\documentclass[12pt,reqno]{amsart}
\usepackage{amsmath,amssymb,amsfonts,amscd,latexsym,amsthm,mathrsfs}
\usepackage{color}
\usepackage{relsize}
\usepackage{hyperref}
\usepackage{fullpage}
\usepackage{enumitem}  

\usepackage{mathtools}
\usepackage{tikz}

\newcommand{\bad}{\textbf{Bad}}   

\newtheorem{theorem}{Theorem}
\newtheorem{proposition}[theorem]{Proposition}        

\newtheorem{remark}[theorem]{Remark}
\newtheorem{lemma}[theorem]{Lemma}
\newtheorem{corollary}[theorem]{Corollary}

\newtheorem{example}[theorem]{Example}
\newtheorem{conjecture}[theorem]{Conjecture} 

\title{Simultaneous approximation to pairs of real numbers}
\subjclass[2020]{11J13,11J70}

\author{Tapani Matala-aho}
\address{Tapani Matala-aho, Mathematics, Aalto University, P.O. Box 11100, FI-00076 Aalto, Finland}
\email{tapani.matala-aho@aalto.fi}

\begin{document}

\begin{abstract}
Let $B_{m}$ and $D_{n}$ be the denominators of the $m$th and $n$th convergent of the real numbers $\alpha$ and $\beta$, respectively.
We introduce an abnormal method in the theory of simultaneous Diophantine approximation to the pair $\alpha,\beta$.
Namely, the question of finding simultaneous approximation is turned into study of small solutions of a linear Diophantine equation
$xB_{m} + yB_{m+1} = zD_{n} + vD_{n+1}$. 
The Thue-Siegel's lemma guarantees the existence of a non-zero integer vector $(x,y,z,v)$ 
in such a way that $|x|,|y|,|z|,|v|$ are bounded above by $\big( B_{m} + B_{m+1} + D_{n} + D_{n+1} \big)^{1/3}$. 
Thereby we can construct an integer 
$q:=xB_{m} + yB_{m+1} = zD_{n} + vD_{n+1}\ge 1$,
a common denominator, which by the theory of continued fractions gives simultaneous approximations to $\alpha$ and $\beta$. 
We give also a variety of explicit constructions without the Thue-Siegel's lemma.
Let $\|\alpha\|:=\underset{k\in\mathbb{Z}}\min\{|\alpha-k|\}$.
For a class of numbers, including particular equivalent numbers $\alpha$ and $\beta$, we show there exist a real number 
$\kappa=\kappa(\alpha,\beta)>1/2$ and infinitely many explicitly constructible positive integers $q$ such that
$\|q\alpha\| \le \frac{1}{q^{\kappa}}$ and $\|q\beta\| \le \frac{1}{q^{\kappa}}$. 
As the result improves Dirichlet's theorem on simultaneous approximation 
it also confirms the classical Littlewood conjecture for such a pair $\alpha,\beta$. 
In addition, we present general criteria for the classical Littlewood conjecture as well as for its $p$-adic counterpart.
\end{abstract}

\maketitle 

\date{}

\section{Introduction}

Dirichtlet's theorem on simultaneous approximation to a system of real numbers is a fundamental piece
in the theory of Diophantine approximation. 
\begin{proposition}\label{Dirichtlet}\cite{SCHMIDT1}
Dirichtlet.
Let at least one of the numbers
$\alpha_1,...,\alpha_m \in\mathbb{R}$
be irrational. 
Then there exist infinitely many integers $q_k\in\mathbb{Z}_{\ge 1}$ and primitive vectors
$\overline{v}_k=(q_k,p_{1,k},...,p_{m,k})^t\in\mathbb{Z}^{m+1}\setminus\{\overline{0}\}$
satisfying
\begin{equation*}\label{linform7}
|q_k\alpha_i+p_{i,k}|<\frac{1}{q_k^{1/m}}
\end{equation*}
for all $i=1,...,m$.
\end{proposition}
Here we are interested in a particular case - simultaneous approximation to pairs of real numbers.
So, let us have two numbers $\alpha,\beta\in\mathbb{R}$, where at least one is irrational.
Then there exist infinitely many $q\in\mathbb{Z}_{\ge 1}$ such that
\begin{equation}\label{Dirit2}
\|q\alpha\| < \frac{1}{q^{1/2}},\quad
\|q\beta\|  < \frac{1}{q^{1/2}}\,,
\end{equation}
where the notation 
$\|\alpha\|:=\underset{k\in\mathbb{Z}}\min\{|\alpha-k|\}$
is devoted to the distance of $\alpha$ from the nearest integer. 
Note, Dirichtlet's theorem is non-effective, it does not give a construction for appropriate integers $q$, 
called common denominators in the sequel.

It is a challenge to improve the upper bound in \eqref{Dirit2} to any pair of real numbers.
However, there exists such pairs. For example, an extremal number $\alpha$ satisfies
\begin{equation*}\label{Extremal1}
 \|q\alpha\|   \le \frac{1}{q^{\tau}},\qquad 	
 \|q\alpha^2\| \le \frac{1}{q^{\tau}},\quad \tau=\frac{\sqrt{5}-1}{2}=0.61803\ldots,
\end{equation*} 
see, Damien Roy \cite{ROY2008}. 
Let $\alpha$ be an algebraic number of degree $k\in\mathbb{Z}_{\ge 1}$ over $\mathbb{Q}$.
The height $H_k(\alpha)$ is defined as the largest absolute value of the coefficients of the minimal polynomial of $\alpha$.
Let $n\in\mathbb{Z}_{\ge 1}$ and let $\xi$ be a transcendental real number.
Denote by $\tau_{\xi}(n)$ the supremum of all real numbers $\tau$ such that $|\xi-\alpha| \le H_k(\alpha)^{-\tau}$
holds for infinitely many algebraic integers $\alpha$ of degree $k\le n$ over $\mathbb{Q}$.
Define $\tau(n)=\sup\{\tau_{\xi}(n)\,|\, \xi\ \text{is a transcendental real number}\}$. 
Roy applied extremal numbers for proving $\tau(3)=\frac{\sqrt{5}+3}{2}=2.61803\ldots$, therefore completing the work 
\cite{DAV_SCH_1969} of Harold Davenport and Wolfgang Schmidt.

On the other hand, Pad\'e approximations of particular generalized hypergeometric series produce simultaneous approximations
rather close to the upper bound in \eqref{Dirit2}. 
For example, the Hermite-Pad\'e approximations of the exponential series imply that
there exist infinitely many explicitly given $q\in\mathbb{Z}_{\ge 1}$ such that 
$\|qe\| \le \frac{1}{q^{\kappa}}$ and $\|qe^2\| \le \frac{1}{q^{\kappa}}$ hold with 
$\kappa=\frac{1}{2}-1/\log\log q$, see Lemmas 7.1. and 7.2. in Ernvall-Hyt\"onen et al. \cite{AMLOTA_2019}.
Despite the results proved with Pad\'e approximations are slightly weaker than Dirichlet's theorem they are in most cases effective
and therefore important for solving Diophantine equations and other Diophantine questions.
See e.g. Michael Bennett \cite{MIKE_1995} and Fel'dman et al. \cite{FEL_NES_1998}.

As shown by Davenport, Schmidt, and Roy possible improvements to Proposition \ref{Dirichtlet} would be valuable not only per se 
but for applications like to the theory of algebraic numbers and in a natural way to the attempts for resolving the Littlewood conjecture.

Let $\alpha,\beta\in\mathbb{R}$. We shall use notations 
$\ell_q(\alpha,\beta) :=\ q\,\|q\alpha\|\,\|q\beta\|$ for $q\in\mathbb{Z}_{\ge 1}$
and 
$\ell(\alpha,\beta) := \underset{q\to\infty}\liminf\, \ell_q(\alpha,\beta) = \underset{q\to\infty}\liminf\, q\,\|q\alpha\|\,\|q\beta\|$. 
The Littlewood conjecture is a long-standing open problem from 1930s which states that 
\begin{conjecture}\label{Littlewood}
For all pairs $\alpha,\beta\in\mathbb{R}$ holds
\begin{equation}\label{Littlewoodconjecture}
\ell(\alpha,\beta) = 0. 
\end{equation} 
\end{conjecture} 
We sometimes call this as the classical Littlewood conjecture in order to emphasize it has a $p$-adic counterpart.
Let $p\in\mathbb{P}$, the set of prime numbers, and $\beta\in\mathbb{R}$. In the $p$-adic case we denote $\ell_{p,q}(\beta) := q\,|q|_p\,\|q\beta\|$
for $q\in\mathbb{Z}_{\ge 1}$ and 
$\ell_{p}(\beta) := \underset{q\to\infty}\liminf\,\ \ell_{p,q}(\beta) = \underset{q\to\infty}\liminf\,\ q\,|q|_p\,\|q\beta\|$.
The $p$-adic Littlewood conjecture states that 
\begin{conjecture}\label{Littlewood}
For all $p\in\mathbb{P}$ and $\beta\in\mathbb{R}$ holds
\begin{equation*}\label{padicLittlewoodconjecture}
\ell_{p}(\beta) = 0.
\end{equation*}
\end{conjecture}

In 2006 Einsiedler, Katok and Lindenstrauss \cite{EKL2006} made a major breakthrough by showing that the Hausdorff dimension
of the set of exceptions for \eqref{Littlewoodconjecture} is zero. 
\begin{proposition}\label{EKLpropintro}\cite{EKL2006}
Let
\begin{equation*}\label{EKL}
\mathcal{E}:= \big\{ (\alpha,\beta)\in\mathbb{R}^2 \big|\, \ell(\alpha,\beta) > 0 \big\}.
\end{equation*} 
Then the Hausdorff dimension $\dim_H \mathcal{E} = 0$. In fact, $\mathcal{E}$ is a countable union of compact sets with box dimension zero.
\end{proposition}
However, despite Proposition \ref{EKLpropintro}, we still do not know whether, for example, 
$(\sqrt{2},\sqrt{3})$ or $(\alpha,\alpha^{-1})$ for an arbitrary $\alpha$ satisfies \eqref{Littlewoodconjecture}. 
It is well-known that if 
$\ell(\alpha,\beta)>0$, then both $\alpha$ and $\beta$ must lie in the set of badly approximable numbers 
\begin{equation*}\label{}
\bad = \big\{\alpha\in\mathbb{R}\ \big|\ \underset{q\in\mathbb{Z}_{\ge 1}}{\inf}\ q\,\|q\alpha\|>0\big\}.
\end{equation*}
Therefore when considering solely the Littlewood conjecture we may restrict investigations to the pairs of badly approximable numbers. 

Pollington and Velani \cite{POLVEL2000} have shown there are limitations for selecting appropriate denominators $q$ 
for proving $\ell(\alpha,\beta)=0$ for a pair $\alpha,\beta$.

\begin{proposition}\label{PV}\cite{POLVEL2000} 
Given $\alpha\in\bad$ and $\lambda\in(0,1)$, there exists a subset $B_{\lambda}(\alpha)$ of $\bad$ with the Hausdorff dimension
$\dim_H B_{\lambda}(\alpha) = \lambda$, such that for any $\beta\in B_{\lambda}(\alpha)$,
\begin{equation*}\label{POLVELintro}
\|B_n\beta\| \ge \delta \quad\forall\ n\in\mathbb{Z}_{\ge 1},
\end{equation*}
where  $\delta=\delta(\alpha,\lambda)>0$ is a constant, and $B_n$ is the denominator of the $n$th convergent of $\alpha$. 
\end{proposition}

Let $B_{m}$ and $D_{n}$ be the denominators of the $m$th and $n$th convergent of the real numbers $\alpha$ and $\beta$, respectively.
Proposition \ref{PV} shows that one can not prove \eqref{Littlewoodconjecture} for all pairs $\alpha,\beta$ just going through
the denominator sequences $(B_{m})$ and $(D_{n})$. 
Therefore something else will be needed -- for example linear combinations.
Just recently in the work \cite{MATOVA2026} the authors introduced linear combinations like
$q_1=rB_m+sB_{m+1}$ and $q_2=tD_n+vD_{n+1}$ with some integers $r,s,t,v$ for the job of possible common denominators $q$. 
It was shown that for all $\alpha,\beta\in\mathbb{R}$ holds
$\underset{q\in\mathbb{Z}_{\ge 1}}\inf\, q\,\|q\alpha\|\,\|q\beta\| < \frac{1}{26}$.
However, the values of the coefficients $r,s,t,v$ were chosen on experimental base with no indication 
how $q_1$ and $q_2$ are related to each other.
As will be explained below, in this work we will marry these linear combinations in order to deliver
a Janus face common denominator for $\alpha$ and $\beta$.

For a more detailed overview and results regarding Diophantine approximations and the Littlewood conjecture, 
we refer to Kristensen \cite{Trends}, Pollington et al. \cite{POLVEL2000} and Schmidt \cite{SCHMIDT1}.

\section{Continued fractions}\label{chapter11}

Most of the following basics on continued fractions may be found from the classic by Oskar Perron \cite{Perron_1913},  
see also \cite{Alfetal_2014} and \cite{HAR_WRI_1989}. For the matrix presentation we refer to \cite{Alfetal_2014}. 

\subsection{Notations}

Every irrational real number $\alpha$ has a unique infinite simple continued fraction representation
\[
\alpha = b_0+\frac{1}{b_1+\dfrac{1}{b_2+_{\ddots}}} = [b_0;b_1,b_2,\ldots],
\]
where $b_0\in\mathbb{Z}$ and $b_n\in\mathbb{Z}_{\ge 1}$ for $n\in\mathbb{Z}_{\ge 1}$. 
The integers $b_n$ are called partial quotients or partial denominators. 
The truncated continued fraction
\[
\frac{A_n}{B_n} := b_0+\frac{1}{b_1+\dfrac{1}{b_2+_{\ddots} \underset{{}+ \dfrac{1}{b_n}}{} }} 
=[b_0;b_1,b_2,\ldots,b_n] 
\]
is called the $n$th convergent of $\alpha$. 
The numerators $A_n$ and the denominators $B_n$ of the convergents are determined via the recurrence formulae
\begin{equation}\label{recurrences}
\begin{cases}
A_{n+2}=b_{n+2}A_{n+1}+A_{n}, \\
B_{n+2}=b_{n+2}B_{n+1}+B_{n},\quad n\in\mathbb{Z}_{\ge 0},
\end{cases}
\end{equation}
with initial values $A_0=b_0$, $B_0=1$, $A_1=b_0b_1+1$ and $B_1=b_1$. Consistent with recurrences \eqref{recurrences}
we may set $A_{-1}=1$, $B_{-1}=0$, $A_{-2}=0$, $B_{-2}=1$.
Further, the numerators and denominators satisfy a determinant formula
\begin{equation}\label{determinant}
A_{n+1}B_{n}-A_{n}B_{n+1}=(-1)^n,\qquad n\in\mathbb{Z}_{\ge 0}.
\end{equation}

Splitting $A_k$. Let us define a sequence $(S_n)$ by the recurrence $S_{n+2}=b_{n+2}S_{n+1}+S_{n}$ for $n\in\mathbb{Z}_{\ge 0}$
with the initial values $S_{0}=0, S_{1}=1$.
By noting $A_{k}=b_{0}B_{k}+S_{k}$ and assuming $b_{0}\ge 1$ we get estimates
\begin{equation}\label{A=B+smallseq}
b_{0}B_{k} < A_{k} < (b_{0}+1)B_{k},\quad k\in\mathbb{Z}_{\ge 1}. 
\end{equation}

Bounds for $A_k$ and $B_k$. By recurrences \eqref{recurrences} it can be shown that
\begin{equation}\label{lowerboundAB}
A_{k} \ge \phi^{k-2},\quad B_{k} \ge \phi^{k-1}
\end{equation}
for $k\ge 2$, where $\phi=\frac{1+\sqrt{5}}{2}$.

The notation
$\alpha_k:=[b_k,b_{k+1},\ldots]$
denotes the $k$th tail of $\alpha=[b_0;b_{1},\ldots]$.
Readily
\begin{equation*}\label{}
\alpha:=\alpha_0=[b_0;b_{1},\ldots]=[b_0;b_{1},\ldots,b_n,\alpha_{n+1}].
\end{equation*}
Note that $[b_0;b_{1},\ldots,b_n,\alpha_{n+1}]$ is not necessarily simple but the recurrences \eqref{recurrences} still hold
for the convergents 
\begin{equation*}\label{}
\frac{\hat A_{k}}{\hat B_{k}} := [b_0;b_{1},\ldots,b_k]=\frac{A_{k}}{B_{k}},\quad k=0,1,\ldots,n,\quad 
\frac{\hat A_{n+1}}{\hat B_{n+1}} := [b_0;b_{1},\ldots,b_n,\alpha_{n+1}].
\end{equation*}
Therefore
\begin{equation*}\label{virhetermi}
\alpha=[b_0;b_{1},\ldots,b_n,\alpha_{n+1}]=\frac{\hat A_{n+1}}{\hat B_{n+1}}=
\frac{\alpha_{n+1}A_n+A_{n-1}}{\alpha_{n+1}B_n+B_{n-1}},
\end{equation*}
which readily implies
\begin{equation*}\label{}
\begin{split}
\alpha - \frac{A_{n}}{B_{n}} & = \frac{(-1)^n}{B_n(\alpha_{n+1}B_n+B_{n-1})}
                               = \frac{(-1)^n}{B_n^2(\alpha_{n+1}+B_{n-1}/B_n)}\\
                             & = \frac{(-1)^n}{B_n^2([b_{n+1},b_{n+2},\ldots]+[0;b_{n},b_{n-1},\ldots,b_{1}])}.
\end{split}
\end{equation*} 
Consequently
\begin{equation*}\label{}
\frac{1}{ B_n((b_{n+1}+1)B_n+B_{n-1})}  \le \left|\alpha - \frac{A_{n}}{B_{n}} \right|  \le \frac{1}{B_n(b_{n+1}B_n+B_{n-1})},
\end{equation*}
where the upper bound may be estimated by
\begin{equation*}\label{favoritebound}
\frac{1}{B_n(b_{n+1}B_n+B_{n-1})} = \frac{1}{B_nB_{n+1}} \le \frac{1}{B_n(B_n+B_{n-1})} \le \frac{1}{B_n^2}.
\end{equation*}

For the nearest integer function holds $\|k\pm\alpha\| = \|\alpha\|$ for every $k\in\mathbb{Z}$. 
In particular, we have
\begin{equation}\label{Bn1perbound}
\|B_n\alpha\| = \|B_n\alpha - A_n\| \le \frac{1}{B_{n+1}}\,,
\end{equation} 
if $|B_n\alpha - A_n|<1/2$.
Moreover
\begin{equation*}\label{lyhinet}
\|\gamma + \delta\| = |\gamma + \delta| \le |\gamma| + |\delta| = \|\gamma\| + \|\delta\|,
\end{equation*}
if $|\gamma|<1/4$, $|\delta|<1/4$. 
Let $x,y\in\mathbb{Z}$ and $q=xB_{m}+yB_{m+1}$. Then we get a crucial tool
\begin{equation}\label{qalphabetabound}
\|q\alpha\|  =  \|xB_{m}\alpha + yB_{m+1}\alpha\| = \|x(B_{m}\alpha-A_{m}) + y(B_{m+1}\alpha-A_{m+1})\| 
            \le \frac{|x|}{B_{m+1}} + \frac{|y|}{B_{m+2}},
\end{equation} 
if $\frac{|x|}{B_{m+1}}<1/4$ and $\frac{|y|}{B_{m+2}}<1/4$.

\subsection{Matrices}

There is one to one correspondence between the convergent 
\begin{equation*}\label{}
\frac{A_n}{B_n} = [r_{0},r_{1},\ldots,r_{n}]							
\end{equation*}
and the matrix product
\begin{equation*}\label{}
\begin{pmatrix}
A_{n} & A_{n-1} \\
B_{n} & B_{n-1}
\end{pmatrix}
= 
\mathcal{R}_{0}\mathcal{R}_{1}\cdots\mathcal{R}_{n},\quad
\mathcal{R}_{k} :=
\begin{pmatrix}
r_{k} & 1 \\
1     & 0
\end{pmatrix}, 
\end{equation*}
for $k=0,1,\ldots,n$. 
In particular,
\begin{equation*}\label{}
\begin{pmatrix}
A_{0} & A_{-1} \\
B_{0} & B_{-1}
\end{pmatrix}
= 
\begin{pmatrix}
r_{0} & 1 \\
1     & 0
\end{pmatrix}\,. 
\end{equation*}
As usual the empty matrix product is defined to be the identity matrix, say $I$.

Here we note that
\begin{equation*}\label{}
\begin{pmatrix}
r_{k} & 1 \\
1     & 0
\end{pmatrix}^{-1}
=
\begin{pmatrix}
0 & 1      \\
1 & - r_{k}
\end{pmatrix}, 
\quad
\begin{pmatrix}
A_{n} & A_{n-1} \\
B_{n} & B_{n-1}
\end{pmatrix}^{-1}
=
(-1)^{n-1}
\begin{pmatrix}
  B_{n-1} & - A_{n-1} \\
  - B_{n} & A_{n}
  \end{pmatrix}\,.
\end{equation*}
Further, write
\begin{equation*}\label{}
\begin{pmatrix}
a_{11} & a_{12}  \\
a_{21} & a_{22} 
\end{pmatrix}
:=
\begin{pmatrix}
a & 1 \\
1 & 0
\end{pmatrix}
\begin{pmatrix}
b & 1 \\
1 & 0
\end{pmatrix}
=
\begin{pmatrix}
ab+1 & a \\
b    & 1
\end{pmatrix}\,.
\end{equation*}
If $a,b\ge 1$, then
$a_{11} \ge a_{12}, a_{21}, a_{22}$.
More generally write
\begin{equation}\label{genprod}
\begin{pmatrix}
a_{11} & a_{12}  \\
a_{21} & a_{22} 
\end{pmatrix}
:= 
\begin{pmatrix}
r_{1} & 1 \\
1     & 0
\end{pmatrix} 
\cdots
\begin{pmatrix}
r_{k} & 1 \\
1     & 0
\end{pmatrix}\,.
\end{equation}
If $r_i\ge 1$ for $i=1,\ldots,k$, then
\begin{equation}\label{11termbiggest}
a_{11} \ge a_{12}, a_{21}, a_{22}. 
\end{equation}
Let
\begin{equation*}\label{}
\mathcal{S} = 
\begin{pmatrix}
a_{11} & a_{12}  \\
a_{21} & a_{22} 
\end{pmatrix},
\quad
\mathcal{T} = 
\begin{pmatrix}
b_{11} & b_{12}  \\
b_{21} & b_{22} 
\end{pmatrix}
\end{equation*}
be matrices of the form \eqref{genprod} and satisfy condition \eqref{11termbiggest}. Then the product
\begin{equation*}\label{}
\mathcal{S}\cdot\mathcal{T} = 
\begin{pmatrix}
c_{11} & c_{12}  \\
c_{21} & c_{22} 
\end{pmatrix}
\end{equation*}
satisfies condition \eqref{11termbiggest}, too. Further,
\begin{equation}\label{11termestimates}
a_{11}b_{11} \le  c_{11} \le 2a_{11}b_{11}.
\end{equation}


\subsection{Periodic type continued fractions}

Let 
\begin{equation*}\label{}
\alpha = \big[ \overline{b_0,b_1,\ldots,b_{J-1}} \big],\quad J\in\mathbb{Z}_{\ge 1},							
\end{equation*}
be a purely periodic continued fraction with a period $J$.
Write
\begin{equation*}\label{ABJakso}
\begin{pmatrix}
A_{J-1} & A_{J-2} \\
B_{J-1} & B_{J-2}
\end{pmatrix}
:=
\mathcal{R}_{0}\cdots\mathcal{R}_{J-1},\quad 
\mathcal{R}_n=
\begin{pmatrix}
b_n  & 1 \\
1    & 0
\end{pmatrix}\,.
\end{equation*}
Then
\begin{equation*}\label{}
\begin{split}
  \begin{pmatrix}
  A_{(h+1)J-1}  & A_{(h+1)J-2} \\
  B_{(h+1)J-1}  & B_{(h+1)J-2}
  \end{pmatrix}
& = \left(\mathcal{R}_{0}\cdots\mathcal{R}_{J-1}\right)^{h+1} 
  = \mathcal{R}_{0}\cdots\mathcal{R}_{J-1}\cdot\left(\mathcal{R}_{0}\cdots\mathcal{R}_{J-1}\right)^{h} \\
& =
  \begin{pmatrix}
  A_{J-1} & A_{J-2} \\
  B_{J-1} & B_{J-2}
  \end{pmatrix}
	\begin{pmatrix}
  A_{hJ-1} & A_{hJ-2} \\
  B_{hJ-1} & B_{hJ-2}
  \end{pmatrix}\,,
	\end{split}
\end{equation*}
which implies
\begin{equation}\label{ABrelation1}
\begin{split}
A_{(h+1)J-1} & = A_{J-1} A_{hJ-1} + A_{J-2} B_{hJ-1}, \\
B_{(h+1)J-1} & = B_{J-1} A_{hJ-1} + B_{J-2} B_{hJ-1}
\end{split}
\end{equation}
for $h\in\mathbb{Z}_{\ge 1}$.
We also have
\begin{equation*}\label{}
  \begin{pmatrix}
  A_{(h+1)J-1}  & A_{(h+1)J-2} \\
  B_{(h+1)J-1}  & B_{(h+1)J-2}
  \end{pmatrix}
 =
  \begin{pmatrix}
  A_{hJ-1} & A_{hJ-2} \\
  B_{hJ-1} & B_{hJ-2}
  \end{pmatrix}
  \begin{pmatrix}
  A_{J-1} & A_{J-2} \\
  B_{J-1} & B_{J-2}
  \end{pmatrix}\,.
\end{equation*}
It follows
\begin{equation}\label{ABrelation2}
\begin{split}
A_{(h+1)J-1} & = A_{J-1} A_{hJ-1} + B_{J-1} A_{hJ-2}, \\
B_{(h+1)J-1} & = A_{J-1} B_{hJ-1} + B_{J-1} B_{hJ-2} 
\end{split}
\end{equation}
for $h\in\mathbb{Z}_{\ge 1}$. 
As a consequence of the first relations in \eqref{ABrelation1} and \eqref{ABrelation2} we get
\begin{equation}\label{ABrelation3}
A_{J-2} B_{hJ-1} = B_{J-1} A_{hJ-2} 
\end{equation}
for $h\in\mathbb{Z}_{\ge 1}$. 

In the following, we mean that a relation is valid for infinitely many integers $h\in\mathbb{Z}_{\ge 1}$,
if there exists an increasing infinite sequence $(h_k)$ satisfying that relation for all $k\in\mathbb{Z}_{\ge 1}$.
Continued fractions satisfying \eqref{ABrelation3} for infinitely many $h\in\mathbb{Z}_{\ge 1}$
will be called periodic type continued fractions.
For example, palindromic continued fractions belong to the class of periodic type continued fractions, see below.
Further, we will call a continued fraction block-periodic if it contains arbitrary long periodic blocks.
Two different periodic blocks do not need to contain identical periods.
A block-periodic continued fraction may be purely periodic or eventually periodic.

\subsection{Palindromic continued fractions}

Let $n\in\mathbb{Z}_{\ge 1}$. The $n$th convergent is a palindrome, if 
\begin{equation*}\label{}
[r_{0},\ldots,r_{n}]= [r_{n},\ldots,r_{0}]
\end{equation*}
or equivalently
\begin{equation}\label{palindromematrix}
  \begin{pmatrix}
  A_{n}  & A_{n-1} \\
  B_{n}  & B_{n-1}
  \end{pmatrix}
 = \mathcal{R}_{0}\cdots\mathcal{R}_{n} 
  = \mathcal{R}_{n}\cdots\mathcal{R}_{0}
 =
  \begin{pmatrix}
  A_{n} & A_{n-1} \\
  B_{n} & B_{n-1}
  \end{pmatrix}^{T},
\end{equation}
which happens exactly when
\begin{equation}\label{palindromerule}
B_{n} = A_{n-1},
\end{equation}
where $n+1$ is length of palindrome. 

The mainstream definition for a palindromic continued fraction says it contains arbitrary long convergents which are palindromes.
However, we will call such continued fractions purely palindromic to emphasize a continued fraction is eventually palindromic, 
if it contains arbitrary long palindromes.
An eventually palindromic continued fraction may be purely palindromic (palindromic).

A purely palindromic continued fraction satisfies 
$B_{h_k-1} = A_{h_k-2}$, see \eqref{palindromerule}, for infinitely many $h_k$, where $h_k$ is length of $k$th palindrome. 
By setting $J=1$, $A_{-1}=1$ and $B_{0}=1$ into \eqref{ABrelation3} we see that the palindromic relation $B_{h_k-1} = A_{h_k-2}$
is a particular case of \eqref{ABrelation3}. Therefore palindromic continued fractions belong to the class of periodic type
continued fractions. 

Several interesting continued fractions may be defined via words in a monoid, say, $M:=\langle a,b,c,\ldots \rangle$.
We use the notation $\overset{\rightarrow} v_h$ for a word $r_{0}r_{1}\ldots r_{h-1}r_{h}\in M$.
Then $\overset{\leftarrow} v_h$ denotes the reverse word $r_{h}r_{h-1}\ldots r_{1}r_{0}\in M$.

In the following $\sigma_F$ and $\sigma_T$ are morphisms in the monoid $\langle a,b \rangle$. 
The Fibonacci word may be defined by a limit $\sigma_F^{\infty}(a)$ of the sequence generated by
\begin{equation*}\label{}
\sigma_F(a):=ab,\quad \sigma(b):=a,\quad \sigma_F^{\infty}(a) = abaababa\ldots.
\end{equation*}
Or equivalently, by a limit $\lim \overset{\rightarrow} w_{k}$ of the sequence
\begin{equation*}\label{}
\overset{\rightarrow} w_1:=ab,\quad \overset{\rightarrow} w_{2}:=aba,\quad
\overset{\rightarrow} w_{k+2}:= \overset{\rightarrow} w_{k+1} \overset{\rightarrow} w_{k},\quad k=1,2,\ldots.
\end{equation*}
Analogously the Thue-Morse word can be defined by a limit $\sigma_T^{\infty}(a)$ of the sequence
\begin{equation*}\label{}
\sigma_T(a):=ab,\quad \sigma(b):=ba,\quad \sigma_T^{\infty}(a) = abbabaab\ldots.
\end{equation*}
Let $\ \widehat{}\ $ be a morphism 
\begin{equation*}\label{}
\widehat{a}:=b,\quad \widehat{b}:=a,
\end{equation*}
in the monoid $\langle a,b \rangle$. 
Then the Thue-Morse word may be defined equivalently by a limit $\lim \overset{\rightarrow} w_{k}$ of the sequence
\begin{equation*}\label{}
\overset{\rightarrow} w_1:=ab,\quad \overset{\rightarrow} w_{k+1}
:= \overset{\rightarrow} w_{k} \widehat{\overset{\rightarrow} w_{k}},\quad k=1,2,\ldots.
\end{equation*}

Let now $a,b\in\mathbb{Z}_{\ge 1}$, $a\ne b$.
So we may define Fibonacci and Thue-Morse continued fractions by setting
\begin{equation*}\label{}
\alpha_F:=[\sigma_F^{\infty}(a)] := [a;b,a,a,b,a,b,a,\ldots],							
\end{equation*}
\begin{equation*}\label{}
\alpha_T:=[\sigma_T^{\infty}(a)] := [a;b,b,a,b,a,a,b,\ldots],							
\end{equation*}
respectively.
It is well-known that the Fibonacci word and the Thue-Morse word are purely palindromic.

If both $\alpha$ and $\beta$ belong to the same real quadratic field, they satisfy the Littlewood conjecture.
The values of eventually periodic continued fractions are quadratic irrationals. 
Further, if $\alpha$ belongs to a real quadratic field, then all numbers equivalent to $\alpha$ belong to the same real quadratic field.
Thus, when considering the Littlewood conjecture alone it would be interesting to know whether $\alpha$ and $\beta$ 
are quadratic irrationals or not. Towards this direction Martine Queff\'elec \cite{QUO_1998} has proved 
the Thue-Morse continued fraction is transcendental.
Further, Boris Adamczewski and Yann Bugeaud \cite{Boris_Yann_2007} proved that palindromic continued fractions are either 
quadratic irrationals or transcendental.
\begin{proposition}\label{}\cite{Boris_Yann_2007}
Let $(a_{\ell})_{\ell\ge 0}$ be a sequence of positive integers. If the word 
$\overset{\rightarrow} a =a_0a_1a_2\ldots$ begins in arbitrarily long palindromes, then the real number 
$\alpha:=[a_0;a_1,a_2,\ldots]$ is either quadratic irrational or transcendental.
\end{proposition}
It is known that the Fibonacci word is non-periodic. Hence, the Fibonacci continued fraction is transcendental, too. 

Throughout the entire work we will identify an irrational real number with its simple continued fraction expansion.
For example, if we say $\alpha$ is palindromic number it means the simple continued fraction expansion of $\alpha$ is palindromic.

\section{A glimpse to explicit constructions}

Now it is a time to introduce our method which allows to improve Dirichtlet's theorem for certain pairs of real numbers
as well as to introduce specific criterions for the Littlewood conjecture. 
Let $B_{m}$ and $D_{n}$ be the denominators of the $m$th and $n$th convergent of the real numbers $\alpha$ and $\beta$,
respectively. Then there exist integers $x,y,z,v$ such that $q:=xB_{m} + yB_{m+1} = zD_{n} + vD_{n+1}\ge 1$
and $|x|,|y|,|z|,|v|$ are bounded above by $\big( B_{m} + B_{m+1} + D_{n} + D_{n+1} \big)^{1/3}$.
A general method is based on Thue-Siegel's lemma but we give also a variety of explicit constructions.

Considering explicit constructions we will pay a particular attention to a class of periodic type continued fractions, 
which contain for example palindromic continued fractions.  
In addition, we construct examples of block-periodic continued fractions as well as eventually palindromic continued fractions. 

Boris Adamczewski and Yann Bugeaud \cite{Boris_Yann_2006} have shown that equivalent palindromic numbers 
$\alpha$ and $\beta$ satisfy  
\begin{equation}\label{BORISYANNbound}
\underset{q\to\infty}\liminf\, q^2\,\|q\alpha\|\,\|q\beta\| < \infty
\end{equation}
and a fortiori the Littlewood conjecture, too. Further,
Adamczewski and Bugeaud \cite{Boris_Yann_2007} proved a bunch of qualitative results for palindromic continued fractions.

Our method differs considerable from the work \cite{Boris_Yann_2006} and allows us to study not only equivalent palindromic numbers 
but a wider class of pairs of real numbers in an explicit quantitative manner.
The notations
\begin{equation*}\label{}
\frac{A_m}{B_m}=[b_0;b_{1},\ldots,b_m],\quad \frac{C_n}{D_n}=[d_0;d_{1},\ldots,d_n],\quad \frac{E_k}{F_k}=[f_0;f_{1},\ldots,f_k],
\end{equation*}
are devoted to the convergents of $\alpha=[b_0;b_{1},\ldots]$, $\beta=[d_0;d_{1},\ldots]$ and $\gamma=[f_0;f_{1},\ldots]$, respectively.
Let us study equivalent numbers $\alpha$ and $\beta$ defined by
\begin{equation*}\label{}
\alpha = [b_0;b_1,\ldots,b_{M},f_0,f_{1},\ldots],\quad \beta = [d_0;d_1,\ldots,d_{N},f_0,f_{1},\ldots].							
\end{equation*}

So, our target is to construct relations like $xB_{m} + yB_{m+1} = zD_{n} + vD_{n+1}$ between 
the denominators $B_{m}$, $B_{m+1}$, $D_{n}$ and $D_{n+1}$.
We call such relation a purely external relation connecting $\alpha$ and $\beta$.
 
Now we will describe how that can be done for particular equivalent continued fractions, 
say $\alpha$ and $\beta = [d_0;d_1,\ldots,d_{N},\alpha]$.   
In Section \ref{secEQCF} we prove 
\begin{equation}\label{EQUIVRELVINOintro}
D_{N+L} = D_{N}A_{L-1} + D_{N-1}B_{L-1}
\end{equation}
for all $L\in\mathbb{Z}_{\ge 1}$. Relation \eqref{EQUIVRELVINOintro} is called a mixed external relation 
because it ties the numerator $A_{L-1}$ with the denominators $B_{L-1}$ and $D_{N+L}$.
 
If there were another relation between $A$s and $B$s, then we may try to eliminate the term $A_{L-1}$ in relation \eqref{EQUIVRELVINOintro}.  
For example, a purely period continued fraction $\alpha$ with period $J\in\mathbb{Z}_{\ge 1}$ satisfies 
\begin{equation}\label{ABrelation4intro} 
A_{J-2} B_{hJ-1} = B_{J-1} A_{hJ-2}, 
\end{equation}
for all $h\in\mathbb{Z}_{\ge 1}$. We call such relation internal relation of $\alpha$.

Combining relations \eqref{EQUIVRELVINOintro} and \eqref{ABrelation4intro} at $L-1=hJ-2$, yields to
\begin{equation*}\label{BDmarried}
B_{J-1} D_{N+hJ-1} = A_{J-2} D_{N} B_{hJ-1} + B_{J-1} D_{N-1} B_{hJ-2}.
\end{equation*}
Therefore we may define infinite number of non-zero common denominators
\begin{equation*}\label{qBDmarried}
q_h := B_{J-1} D_{N+hJ-1} = A_{J-2} D_{N} B_{hJ-1} + B_{J-1} D_{N-1} B_{hJ-2},
\end{equation*}
where the coefficients of $D_{N+hJ-1}$, $B_{hJ-1}$ and $B_{hJ-2}$ are constant for fixed $J$ and $N$ while $h$ is a free variable.

For the following theorem we will apply the above like deduction in a little bit more sophisticated manner. 
Namely, see Section \ref{pertype}, relation \eqref{PERTYPEEQUIV}, which produces us a Janus face denominator 
$q_h$ for equivalent periodic type continued fractions $\alpha$ and $\beta$.

\begin{theorem}\label{alphabetagamma}
Let $(h_k)$ be an increasing infinite sequence of positive integers, $J\in\mathbb{Z}_{\ge 1}$, and let  
$\gamma$ be a continued fraction satisfying the relation 
\begin{equation}\label{EFgammarel}
E_{J-2} F_{h_kJ-1} = F_{J-1} E_{h_kJ-2} 
\end{equation}
for all $k\in\mathbb{Z}_{\ge 1}$.
Suppose $\alpha$ and $\beta$ are equivalent continued fractions in the manner 
$\alpha = [b_0;b_1,\ldots,b_M,\gamma]$, $\beta = [d_0;d_1e,\ldots,d_{N},\gamma]$.
Write
$c_{1}=c_{1}(\beta,\gamma):= (F_{J-1} D_{N-1} + E_{J-2} D_{N})^2$, 
$c_{2}=c_{2}(\alpha,\gamma):=(F_{J-1} B_{M-1} + E_{J-2} B_{M})^2$
and
\begin{equation}\label{qhkBDEF}
q_{h_k} := F_{J-1} B_{M-1} D_{N+h_kJ-1} + E_{J-2} B_{M} D_{N+h_kJ} = F_{J-1} D_{N-1} B_{M+h_kJ-1} + E_{J-2} D_{N} B_{M+h_kJ}
\end{equation}
for $k\in\mathbb{Z}_{\ge 1}$.
Then the inequalities
\begin{equation}\label{eqcfbound}
\begin{split}
 \|q_{h_k}\alpha\| & \le \frac{c_{1}}{q_{h_k}}\,,\\ 
 \|q_{h_k}\beta\|  & \le \frac{c_{2}}{q_{h_k}}
\end{split}
\end{equation} 
hold simultaneously for all $k\in\mathbb{Z}_{\ge 1}$.
\end{theorem}

Proof. The construction of $q_{h_k}$ will be delivered in Section \ref{pertype}, see relation \eqref{PERTYPEEQUIV}.
First we give estimates for $q_{h_k}$ in terms of $B_{M+h_kJ}$ and $D_{N+h_kJ}$ as follows
\begin{equation*}\label{qBDbounds}
\begin{split}
q_{h_k} & = F_{J-1} D_{N-1} B_{M+h_kJ-1} + E_{J-2} D_{N} B_{M+h_kJ} \le (F_{J-1} D_{N-1} + E_{J-2} D_{N}) B_{M+h_kJ}, \\
q_{h_k} & = F_{J-1} B_{M-1} D_{N+h_kJ-1} + E_{J-2} B_{M} D_{N+h_kJ} \le (F_{J-1} B_{M-1} + E_{J-2} B_{M}) D_{N+h_kJ}.
\end{split}
\end{equation*}
On the other hand, by the theory of continued fractions, see Section \ref{chapter11} equation \eqref{qalphabetabound}, we get the bounds
\begin{equation*}\label{}
\|q\alpha\| \le \frac{|x|}{B_{m+1}} + \frac{|y|}{B_{m+2}}\,,\quad
\|q\beta\|  \le \frac{|z|}{D_{n+1}} + \frac{|v|}{D_{n+2}}
\end{equation*} 
for $q=xB_{m} + yB_{m+1}=zD_{n} + vD_{n+1}$.
Thereby
\begin{equation*}\label{}
\begin{split}
 \|q_{h_k}\alpha\| & =  \| F_{J-1} D_{N-1} B_{M+h_kJ-1}\,\alpha + E_{J-2} D_{N} B_{M+h_kJ}\,\alpha\| \\
               & \le \frac{F_{J-1} D_{N-1}}{B_{M+h_kJ}} + \frac{E_{J-2} D_{N}}{ B_{M+h_kJ+1}} 
                 \le \frac{F_{J-1} D_{N-1} + E_{J-2} D_{N}}{B_{M+h_kJ}} 	
							   \le \frac{(F_{J-1} D_{N-1} + E_{J-2} D_{N})^2}{q_{h_k}}, \\
 \|q_{h_k}\beta\|  & =  \| F_{J-1} B_{M-1} D_{N+h_kJ-1}\,\beta + E_{J-2} B_{M} D_{N+h_kJ}\,\beta \| \\
               & \le \frac{F_{J-1} B_{M-1}}{D_{N+h_kJ}} + \frac{E_{J-2} B_{M}}{D_{N+h_kJ+1} }
                 \le \frac{F_{J-1} B_{M-1} + E_{J-2} B_{M}}{D_{N+h_kJ}}
							   \le \frac{(F_{J-1} B_{M-1} + E_{J-2} B_{M})^2}{q_{h_k}}.\qed								
\end{split}
\end{equation*}

Here we note that all the common denominators $q_{h_k}$, $k\in\mathbb{Z}_{\ge 1}$, in \eqref{qhkBDEF} are totally explicit.
As $J,M$ and $N$ are fixed, then $c_1$ and $c_2$ are constant with respect to $k$.
Therefore inequalities \eqref{eqcfbound} surpass Dirichtlet's theorem, Proposition \ref{Dirichtlet}, 
on simultaneous approximation to the pair $\alpha$, $\beta$.
Note also, that the bounds are totally explicit.

\begin{corollary}\label{alphabetagammacor}
Let $\alpha = [b_0;b_1,\ldots,b_M,\gamma]$, $\beta = [d_0;d_1,\ldots,d_{N},\gamma]$ satisfy 
the assumptions of Theorem \ref{alphabetagamma}. Then
\begin{equation*}\label{}
\ell_{q_{h_k}}(\alpha,\beta) =  q_{h_k}\,\|q_{h_k}\alpha\|\,\|q_{h_k}\beta\| \le \frac{c_{1}c_{2}}{q_{h_k}}
\end{equation*} 
for all $q_{h_k}$. Therefore $\alpha$ and $\beta$ satisfy the Littlewood conjecture $\ell(\alpha,\beta) = 0$.
\end{corollary}

Consequently we receive a totally explicit bound
\begin{equation*}\label{}
\underset{q\to\infty}\liminf\, q^2\,\|q\alpha\|\,\|q\beta\| \le c_{1}(\beta,\gamma)c_{2}(\alpha,\gamma),
\end{equation*} 
where $\gamma$ satisfies relation \eqref{EFgammarel}. 
Adamczewskis and Bugeaud's non-explicit bound \eqref{BORISYANNbound} holds for purely palindromic 
continued fractions $\gamma$ (in our notation), while the explicit bounds in
Theorem \ref{alphabetagamma} and Corollary \ref{alphabetagammacor} are valid not only for periodic and 
purely palindromic continued fractions $\gamma$ but for all continued fractions $\gamma$ which satisfy relation \eqref{EFgammarel}.

To be more precise we now assume  $\gamma$ is a purely palindromic continued fraction with  
an increasing infinite sequence of palindromes of length $h_k$. 
By noting that the palindromic relation
\begin{equation*}\label{}
F_{h_k-1} = E_{h_k-2} 
\end{equation*}
is equal to \eqref{EFgammarel} with $J=1$, $E_{-1}=1$ and $F_{0}=1$ we obtain the following corollaries. 

\begin{corollary}\label{}
Let $\gamma$ be a purely palindromic continued fraction and let $\alpha = [b_0;b_1,\ldots,b_M,\gamma]$ and 
$\beta = [d_0;d_1,\ldots,d_{N},\gamma]$ satisfy the assumptions of Theorem \ref{alphabetagamma}. Then
\begin{equation*}\label{eqcfbound3}
\begin{split}
 \|q_{h_k}\alpha\| & \le \frac{(D_{N-1} + D_{N})^2}{q_{h_k}}\,,\\ 
 \|q_{h_k}\beta\|  & \le \frac{(B_{M-1} + B_{M})^2}{q_{h_k}}
\end{split}
\end{equation*} 
hold simultaneously for all 
\begin{equation*}\label{qhkBDEF3}
q_{h_k} := B_{M-1} D_{N+h_k-1} + B_{M} D_{N+h_k} = D_{N-1} B_{M+h_k-1} + D_{N} B_{M+h_k},\quad k\in\mathbb{Z}_{\ge 1},
\end{equation*}
where $h_k$ is the length of the $k$th palindrome in $\gamma$. 
\end{corollary}  

Next we consider a special case, where $\alpha=\gamma$, $M=-1$ and $\beta=\alpha^{-1}=[0;\gamma]$, $N=0$.
That forces $B_{M-1}=B_{-2}=1$, $B_{M}=B_{-1}=0$, $D_{N-1}=D_{-1}=0$ and $D_{N}=D_{0}=1$. 

\begin{corollary}\label{}
Let $\alpha$ be a palindromic continued fraction.
Then 
\begin{equation*}\label{eqcfbound4}
\|q_{h_k}\alpha\| \le \frac{1}{q_{h_k}}\,,\quad 
\|q_{h_k}\alpha^{-1}\| \le \frac{1}{q_{h_k}}
\end{equation*} 
hold simultaneously and
\begin{equation*}\label{} 
\ell_{q_{h_k}}(\alpha,\alpha^{-1}) = \ q_{h_k}\,\|q_{h_k}\alpha\|\,\|q_{h_k}\alpha^{-1}\| \le \frac{1}{q_{h_k}}
\end{equation*}
for all 
$q_{h_k} := D_{N+h_k-1} = B_{M+h_k},\ k\in\mathbb{Z}_{\ge 1}$.
Consequently
\begin{equation*}\label{}
\underset{q\to\infty}\liminf\, q^2\,\|q\alpha\|\,\|q\alpha^{-1}\| \le 1.
\end{equation*} 
\end{corollary} 

Of course, we may study equivalent continued fractions 
$\alpha = [b_0;b_1,\ldots,b_M,\gamma]$, $\beta = [d_0;d_1,\ldots,d_{N},\gamma]$,
where $\gamma$ is any continued fraction satisfying a relation, say,
$xE_{m} + yE_{m+1} = zF_{n} + vF_{n+1}$ instead of \eqref{EFgammarel}.
Therefore $\gamma$ do not necessarily be periodic type or palindromic.
If $x,y,z,v$ are constant wrt $m,n$, then we will get simultaneous upper bounds as in \eqref{eqcfbound} with different constants. 
If $x,y,z,v$ are not constant wrt $m,n$, let $|x|,|y|,|z|,|v|$ be bounded above by 
$\big( B_{m} + B_{m+1} + D_{n} + D_{n+1} \big)^{\tau}$, where $0\le \tau<1/3$.
In this case we still would beat the upper bound $q^{-1/2}$ in Dirichlet's theorem.
Particularly, this works for a class of equivalent block-periodic continued fractions as well as 
for eventually palindromic continued fractions, see sections \ref{blockperiodex} and \ref{eventuallypalindex}.
For both cases we shall construct descriptive examples such that there exist a real number 
$\kappa=\kappa(\alpha,\beta) >1/2$ and infinitely many positive integers $q$ such that
$\|q\alpha\| \le \frac{1}{q^{\kappa}}$ and $\|q\beta\| \le \frac{1}{q^{\kappa}}$.
Evidently the Littlewood conjecture holds for such a pair $\alpha$, $\beta$.

\section{Thue-Siegel's lemma and linear combinations}\label{Chaptersiegel}

The celebrated Thue-Siegel's lemma or Siegel's lemma gives an estimate for the size of an integer solution  
to a system of homogeneous integer equations.
The following version is along Kurt Mahler \cite{MAHLER1975}.

\begin{proposition}\label{SiegelMahler} 
Let $M,N\in\mathbb{Z}_{\ge 1}$, $M<N$ and $a_{mn}\in \mathbb{Z}$ for $m=1,...,M$ and $n=1,...,N$.
Assume that
$V_m:=\sum_{n=1}^{N}|a_{mn}|\ge 1$ for $m=1,...,M$. Then  the system of equations
\begin{equation*}\label{SIEGELEQUATIONS}
\begin{cases}
&a_{11}x_1+a_{12}x_2+...+a_{1N}x_N=0,\\
&a_{21}x_1+a_{22}x_2+...+a_{2N}x_N=0,\\
&...\\
&a_{M1}x_1+a_{M2}x_2+...+a_{MN}x_N=0,
\end{cases}
\end{equation*}
has a non-zero integer solution 
$\overline{z}=(z_1,...,z_N)^t\in\mathbb{Z}^{N}\setminus\{\overline{0}\}$ with
\begin{equation*}\label{SIEGELESTIMATE}
1\le \underset{1\le n\le N}\max|z_n|\le V:=\left\lfloor (V_1\cdots V_M)^{\frac{1}{N-M}} \right\rfloor.
\end{equation*}
\end{proposition}

Of course, Thue-Siegel's lemma shows only the existence of a solution.
However, a solution may be found on experimental base, say, just running $z_1,\ldots,z_N$ through the interval $[-V,V]$. 
Hence, such a solution is effectively determined. 

In the sequel the notation $\|\overline{S}\|_{\infty} = \underset{i=1,...,k}\max|x_{i}|$ will be used for the maximum norm of
$\overline{S} = (x_{1},\ldots,x_{k})$.

\subsection{The classical case}\label{Chapterclassical}

By Proposition \ref{SiegelMahler} there exists a non-zero effectively determined integer solution 
$\overline{S}=(x,y,z,v)^t\in\mathbb{Z}^{4}\setminus\{\overline{0}\}$
of the equation
\begin{equation}\label{pohjayhtalo}
xB_{m} + yB_{m+1} = zD_{n} + vD_{n+1}
\end{equation}
such that
\begin{equation*}\label{BD13}
\|\overline{S}\|_{\infty} 
\le J = J(m,n) := \left\lfloor \big( B_{m} + B_{m+1} + D_{n} + D_{n+1} \big)^{1/3} \right\rfloor
\le \big(4\max\{B_{m+1},D_{n+1}\}\big)^{1/3}.
\end{equation*}

\begin{lemma}\label{Lemmaorxyzv} 
Let $\overline{S}(m,n)=(x,y,z,v)^t$ be a solution of equation \eqref{pohjayhtalo} such that
\begin{equation*}\label{xyzvleBD}
\begin{split}
&1\le |x| < B_{m+1}\quad\text{or}\quad 1\le |y| < B_{m}
\quad\text{or}\quad
1\le |z| < D_{n+1}\quad\text{or}\quad 1\le |v| < D_{n}.
\end{split}
\end{equation*} 
Then
$xB_{m} + yB_{m+1} = zD_{n} + vD_{n+1} \ne 0$,
and we may assume wlog that 
\begin{equation*}\label{qge1}
q := xB_{m} + yB_{m+1} = zD_{n} + vD_{n+1} \ge 1.
\end{equation*}
\end{lemma}\label{}

Proof. By the theory of continued fractions $B_m\perp B_{m+1}$. 
Assume for example that $1\le |x| < B_{m+1}$.
If $q=0$, then
\begin{equation*}\label{Bonnolla}
B_{m}x + B_{m+1}y = 0 \quad\Rightarrow\quad x=kB_{m+1},\ y=-kB_{m},\ k\in\mathbb{Z}.
\end{equation*}
A contradiction.\qed

\begin{lemma}\label{LemmahJD}
Let $\overline{S}(m,n)=(x,y,z,v)^t$ be a non-zero solution of equation \eqref{pohjayhtalo} such that
\begin{equation*}\label{xleDn1}
\|\overline{S}(m,n)\|_{\infty} \le H J(m,n),
\end{equation*}
where $H\ge 1$ is a constant.
If
\begin{equation}\label{happens2}
\begin{split}
B_{m+1} < D_{n},\quad 4^4H^3D_{n+1} <  B_{m}^3,
\end{split}
\end{equation}
then
\begin{equation}\label{q2ge1}
q = xB_{m} + yB_{m+1} = zD_{n} + vD_{n+1} \ge 1.
\end{equation}
Further
\begin{equation*}\label{xB14}
\frac{|x|}{B_{m+1}}<1/4,\quad \frac{|y|}{B_{m+2}}<1/4.
\end{equation*}
\end{lemma}

Proof. By \eqref{happens2} we get
$
H \left(4\max\{B_{m+1},D_{n+1}\}\right)^{1/3} < B_{m}/4.
$
Thus
\begin{equation*}\label{xleDn1}
 \max\{|x|,|y|,|z|,|v|\} \le H J(m,n) \le H \cdot\left(4\max\{B_{m+1},D_{n+1}\}\right)^{1/3} < B_{m}/4.
\end{equation*}
Therefore, Lemma \ref{Lemmaorxyzv} implies \eqref{q2ge1}.
In addition, we conclude $\frac{|x|}{B_{m+1}}<1/4$ and $\frac{|y|}{B_{m+2}}<1/4$.\qed

\subsection{The $p$-adic case}\label{SubChapterpadic}

By Proposition \ref{SiegelMahler} there exists a non-zero effectively determined integer solution 
$\overline{S}_1=(x,z,v)^t\in\mathbb{Z}^{3}\setminus\{\overline{0}\}$
of equation 
\begin{equation}\label{pohjayhtalop}
xp^{m} = zD_{n} + vD_{n+1}
\end{equation}
such that
\begin{equation*}\label{BD13p}
\begin{split}
\|\overline{S}(m,n)\|_{\infty} 
\le J_p = J_p(m,n) :=  \left\lfloor \big( p^{m} + D_{n} + D_{n+1} \big)^{1/2} \right\rfloor 
\le \big( 3 \max\{ p^{m}, D_{n+1} \} \big)^{1/2}.
\end{split}
\end{equation*}

\begin{lemma}\label{Lemmaorxzv} 
Let $\overline{S}(m,n)=(x,z,v)^t$ be a solution of equation \eqref{pohjayhtalop} such that
\begin{equation*}\label{xyzvleBDp}
1\le |z| < D_{n+1}\quad\text{or}\quad 1\le |v| < D_{n}.
\end{equation*} 
Then
$xp^{m} = zD_{n} + vD_{n+1} \ne 0$,
and we may assume wlog that 
\begin{equation*}\label{qge1p}
q := xp^{m} = zD_{n} + vD_{n+1} \ge 1.
\end{equation*}
\end{lemma}\label{}

Proof. Similar to the proof of Lemma \ref{Lemmaorxyzv}.

\begin{lemma}\label{LemmahJDp}
Let $\overline{S}(m,n)=(x,z,v)^t$ be a non-zero solution of equation \eqref{pohjayhtalop} such that
\begin{equation*}\label{xleDn1p}
\|\overline{S}(m,n)\|_{\infty} \le H J_p(m,n),
\end{equation*}
where $H\ge 1$ is a constant.
If
\begin{equation*}\label{happens2p}
p^{m} < D_{n},\quad 
3\cdot 4^2H^2 D_{n+1} < D_{n}^2,
\end{equation*}
then
$q = xp^{m} = zD_{n} + vD_{n+1} \ge 1$.
In addition
\begin{equation*}\label{xB14p}
\frac{|z|}{D_{n+1}}<1/4,\quad \frac{|v|}{D_{n+2}}<1/4.
\end{equation*}
\end{lemma}

Proof. Similar to the proof of Lemma \ref{LemmahJD}.

\section{Criterions for the Littlewood conjecture}

\subsection{The classical case}

Let 
\begin{equation*}\label{}
\frac{A_m}{B_m}=[b_0;b_{1},\ldots,b_m],\quad \frac{C_n}{D_n}=[d_0;d_{1},\ldots,d_n]
\end{equation*}
be convergents of $\alpha=[b_0;b_{1},\ldots]$ and $\beta=[d_0;d_{1},\ldots]$, respectively.
First we give some preliminaries for studying the quantity $\ell(\alpha,\beta)$.
For $(m,n)\in\mathbb{Z}^{2}_{\ge 1}$ let 
$\overline{S}(m,n)=(x,y,z,v)=(x_{m,n},y_{m,n},z_{m,n},v_{m,n})\in\mathbb{Z}^{4}\setminus\{\overline{0}\}$
be a non-zero integer solution of the equation
\begin{equation}\label{pohjayhtalo1}
xB_{m} + yB_{m+1} = zD_{n} + vD_{n+1}.
\end{equation}
By the Thue-Siegel's lemma, see Chapter \ref{Chaptersiegel}, we find a non-zero integer solution
$\overline{S}(m,n)\in\mathbb{Z}^{4}\setminus\{\overline{0}\}$ for equation \eqref{pohjayhtalo1} such that
\begin{equation*}\label{xyzvbounded}
\left\| \overline{S}(m,n) \right\|_{\infty} = \max\{|x|,|y|,|z|,|v|\} 
\le \left\lfloor \big( B_{m} + B_{m+1} + D_{n} + D_{n+1} \big)^{1/3} \right\rfloor =: J(m,n). 
\end{equation*}

For a given function $a:\mathbb{Z}_{\ge 1}\to \mathbb{Z}_{\ge 1}$ we shall consider such solutions
$\overline{S}(a(n),n)=$ $(x_{a(n),n},y_{a(n),n},z_{a(n),n},v_{a(n),n})$ of the equation
\begin{equation}\label{pohjayhtaloank}
xB_{a(n)} + yB_{a(n)+1} = zD_{n} + vD_{n+1}
\end{equation}
that there exists an increasing infinite sequence $(n_k)_{k=1}^{\infty}$ satisfying
\begin{equation*}\label{qkclassical}
q_k := x_{a(n_k),n_k} B_{a(n_k)} + y_{a(n_k),n_k} B_{a(n_k)+1} = z_{a(n_k),n_k} D_{n_k} + v_{a(n_k),n_k} D_{n_k+1} \ge 1
\end{equation*}
for all $k\in\mathbb{Z}_{\ge 1}$. 
Solutions satisfying $z+v\ge 1$, $z,v\ge 0$ will be denoted by $\overline{S}^{++}(m,n)=(x,y,z,v)$ and called ++-solutions.

\begin{theorem}\label{them2}
Let $a:\mathbb{Z}_{\ge 1}\to \mathbb{Z}_{\ge 1}$ be a function such that
\begin{equation}\label{anlim0}
\underset{n\to\infty}\lim\, \frac{ B_{a(n)+1} }{ D_{n} } = 0 
\end{equation} 
and the conditions 
\begin{equation}\label{44H3D}
B_{a(n)+1} \le D_{n},\quad 4^4H^3D_{n+1} <  B_{a(n)}^3
\end{equation} 
hold for all $n\ge n_0\ge 1$. 
Assume that there exists an infinite sequence 
$\big( \overline{S}^{++}(a(n_k),n_k) \big)_{k=1}^{\infty}$ 
of ++-solutions of equation \eqref{pohjayhtaloank} such that
\begin{equation}\label{n++aarettomanmonta}
1\le \left\| \overline{S}^{++}(a(n_k),n_k) \right\|_{\infty} \le H\cdot J(a(n_k),n_k)
\end{equation}
for all $k\in\mathbb{Z}_{\ge 1}$ with some absolute constant $H\ge 1$. Then
$\alpha$ and $\beta$ satisfy the Littlewood conjecture $\ell(\alpha,\beta) = 0$.
\end{theorem}

Proof of Theorem \ref{them2}.
First we note Lemma \ref{LemmahJD} with \eqref{44H3D} and\eqref{n++aarettomanmonta} guarantees 
\begin{equation*}\label{qkconstruction2}
q_k := x_{a(n_k),n_k} B_{a(n_k)} + y_{a(n_k),n_k} B_{a(n_k)+1} = z_{a(n_k),n_k} D_{n_k} + v_{a(n_k),n_k} D_{n_k+1} \ge 1
\end{equation*}
for all $n_k\ge n_0$. 
Write now $m_k=a(n_k)$, $q_k = xB_{m_k} + yB_{m_k+1} = zD_{n_k} + vD_{n_k+1}$ and $J=J(a(n_k),n_k)$ for short.
By the definition of $\overline{S}^{++}(m_k,n_k)$ we have $z+v\ge 1$, $z,v\ge 0$. Therefore
\begin{equation}\label{KAIKENTAKANA}
(z + v)D_{n_k} \le zD_{n_k} + vD_{n_k+1} = xB_{m_k} + yB_{m_k+1} \le 2HJB_{m_k+1}.
\end{equation}
On the other hand 
$q_k = xB_{m_k} + yB_{m_k+1} \le 2HJB_{m_k+1}$.
By \eqref{qalphabetabound} and Lemma \ref{LemmahJD} we get
\begin{equation*}\label{}
\|q_k\alpha\| \le \frac{|x|}{B_{m_k+1}} + \frac{|y|}{B_{m_k+2}}\,,\qquad
\|q_k\beta\|  \le \frac{|z|}{D_{n_k+1}} + \frac{|v|}{D_{n_k+2}}\,.	
\end{equation*} 
Thereby the quantity 
$\ell_{q_k}(\alpha,\beta)=q_k\,\|q_k\alpha\|\,\|q_k\beta\|$ 
may be estimated as follows
\begin{equation*}\label{}
\begin{split}
\ell_{q_k}(\alpha,\beta) 
& \le q_k\,\left(\frac{|x|}{B_{m_k+1}} + \frac{|y|}{B_{m_k+2}}\right)\,\left(\frac{|z|}{D_{n_k+1}} + \frac{|v|}{D_{n_k+2}}\right) \\
& \le 2H JB_{m_k+1}\,\left(\frac{HJ}{B_{m_k+1}} + \frac{HJ}{B_{m_k+2}}\right)\,\left(\frac{z}{D_{n_k+1}} + \frac{v}{D_{n_k+2}}\right) \\
& \le 4H^2J^2\,\frac{z+v}{D_{n_k+1}} 
  \le 4H^2J^2\,\frac{2HJB_{m_k+1}}{D_{n_k}D_{n_k+1}}
  =   \frac{8H^3J^3B_{m_k+1}}{D_{n_k}D_{n_k+1}}.   
\end{split}
\end{equation*} 
We also have $J^3\le 4D_{n_k+1}$. Hence,
\begin{equation*}\label{}
\ell_{q_k}(\alpha,\beta) \le \frac{8H^3J^3B_{m_k+1}}{D_{n_k}D_{n_k+1}} 
\le \frac{2^5H^3B_{m_k+1}D_{n_k+1}}{D_{n_k}D_{n_k+1}}
\le 2^5H^3\frac{B_{a(n_k)+1}}{D_{n_k}}
\end{equation*}
holds for all $n\ge n_0$. Therefore
\begin{equation*}\label{}
\ell(\alpha,\beta) = \underset{q\to\infty}\liminf\, \ell_q(\alpha,\beta) 
\le \underset{k\to\infty}\liminf\, \ell_{q_k}(\alpha,\beta) 
\le \underset{k\to\infty}\liminf\, 2^5H^3\frac{B_{a(n_k)+1}}{D_{n_k}} = 0.\qed
\end{equation*}


Note that for any pair $(m,n)$ equation \eqref{pohjayhtalo1} has at least one solution satisfying
\begin{equation*}\label{firstminimaratkaisu}
1\le \left\| \overline{S}(m,n) \right\|_{\infty} \le J(m,n).
\end{equation*}
Actually, there are uncountable many sequences $\left(\overline{S}(a(n),n)\right)_{n=1}^{\infty}$ of such solutions
with varying functions $a:\mathbb{Z}_{\ge 1}\to \mathbb{Z}_{\ge 1}$.
However, we do not know whether there exists an infinite sequence 
$\left(\overline{S}^{++}(a(n_k),n_k)\right)_{n=1}^{\infty}$ 
of type $++$-solutions satisfying \eqref{n++aarettomanmonta}.
This would not be an issue, if we had an estimate like \eqref{KAIKENTAKANA} for 
solutions, where either $z$ or $v$ is negative.

\begin{remark}\label{remarkthebigger}
If $B_{a(n)+1}$ is much smaller than $D_{n}$, there might be no ++-solutions satisfying \eqref{n++aarettomanmonta}.
Namely, if $2B_{a(n)+1} < D_{n}^{2/3}$, then there are no solutions satisfying $z+v\ge 1$, $z,v\ge 0$ and 
\eqref{n++aarettomanmonta} due to
\begin{equation*}\label{}
xB_{a(n)} + yB_{a(n)+1} < D_{n}^{1/3}D_{n}^{2/3} < zD_{n} + vD_{n+1}.
\end{equation*}
Thus, the bigger $B_{a(n)+1}/D_{n}$ is in the equation
\begin{equation*}\label{pohjayhtalo3}
xB_{a(n)} + yB_{a(n)+1} = zD_{n} + vD_{n+1}
\end{equation*}
the bigger will be the chance for $++$-solutions. 
Let $\kappa:\mathbb{Z}_{\ge 1}\to \mathbb{R}_{>0}$ be a function satisfying the conditions
\begin{equation}\label{kappacond}
\kappa(n)\cdot \log  D_{n+1}\ \underset{n\ \rightarrow\ \infty}\longrightarrow\ \infty,\quad
\kappa(n) \underset{n\ \rightarrow\ \infty}\longrightarrow\ 0.
\end{equation}
Then $B_{a(n)+1}/D_{n}\,\rightarrow\,0$, if 
\begin{equation*}\label{}
B_{a(n)+1} \le D_{n}^{1-\kappa(n)}
\end{equation*}
for all $n$ big enough.
For example, a good choice would be a function like
\begin{equation*}\label{}
\kappa(n) = \frac{\log\log D_{n+1}}{\log D_{n+1}}\,.
\end{equation*}
Of course, there are neverending smaller choices for $\kappa(n)$ satisfying condition \eqref{kappacond}. 
\end{remark} 

As an illustration of condition \eqref{anlim0} we consider an example with the numbers
$\alpha = [0;\overline{1}] = \frac{\sqrt{5}-1}{2}$, $\beta = [0;\overline{2}] = \sqrt{2}-1$.							
Now
$B_1 = 1, B_2 = 2, B_k = B_{k-1} + B_{k-2}$ and 
$D_1 = 2, D_2 = 5, D_k = 2 D_{k-1} + D_{k-2}$,
which may be represented by Binet formulae
\begin{equation*}\label{EFGH}
\begin{split}
B_m &= \frac{1}{\sqrt{5}} \left( \left(\frac{1+\sqrt{5}}{2}\right)^{m+1} +  \left(\frac{1-\sqrt{5}}{2}\right)^{m+1} \right)
    =: \frac{1}{\sqrt{5}} \left( E^{n+1} + F^{n+1} \right),\\
D_n &=  \frac{1}{2\sqrt{2}} \left( \left(1+\sqrt{2}\right)^{n+1} +  \left(1-\sqrt{2}\right)^{n+1} \right)
    =: \frac{1}{2\sqrt{2}} \left( G^{n+1} + H^{n+1} \right)
\end{split}
\end{equation*}
for $n\in\mathbb{Z}_{\ge 1}$.
Denote
\begin{equation*}\label{}
\hat a := \left(1-\kappa(n)\right) \frac{\log G}{\log E} 
=  \left(1-\kappa(n)\right) \frac{\log(1+\sqrt{2})}{\log((1+\sqrt{5})/2)}. 
\end{equation*}
Then
\begin{equation*}\label{}
\underset{n\to\infty} \lim\, \frac{B_{a(n)+1}}{D_n} = 0
\end{equation*}
for all functions $a:\mathbb{Z}_{\ge 1}\to \mathbb{Z}_{\ge 1}$ such that 
$\frac{\hat{a}}{2}n \le a(n)\le \hat{a}n$ for all $n\in\mathbb{Z}_{\ge 1}$. 
In addition, conditions \eqref{44H3D} hold for all $n$ big enough.
Hence, we obtain

\begin{example}\label{exsqrt25}
Let $a:\mathbb{Z}_{\ge 1}\to \mathbb{Z}_{\ge 1}$ be a function such that 
$\frac{\hat{a}}{2}n \le a(n)\le \hat{a}n$ for all $n\in\mathbb{Z}_{\ge 1}$. 
Suppose there exists an infinite sequence of ++-solutions such that
\begin{equation*}\label{n+aarettomanmontaex}
\left\| \overline{S}^{++}(a(n_k),n_k) \right\|_{\infty} \le H J(a(n_k),n_k)
\end{equation*}
for some absolute constant $H\ge 1$. Then
\begin{equation*}\label{} 
\underset{q\to\infty} \liminf\, \ q\,\left\|q(\sqrt{5}+1)/2\right\|\,\left\|q\sqrt{2}\right\| = 0.  
\end{equation*}
\end{example}\label{}

The next criterion does not use assumption \eqref{anlim0}.
Further, we do not need to restrict solutions to be of the type $\overline{S}^{++}(a(n_k),n_k)$.
Instead we assume a slight tightening to the upper bound of $|x|,|y|,|z|,|v|$. 
Also here we assume $B_{a(n)}<D_{n}$ for big enough $n$, which gives a bound $J(a(n),n)^3 \le 4 D_{n+1}$. 
In Theorem \ref{them3} we will apply a modified version of this bound. 
Let again $\kappa:\mathbb{Z}_{\ge 1}\to \mathbb{R}_{>0}$ be a function satisfying the conditions in \eqref{kappacond}.
In the following $J_1(m,n), J_2(m,n)\,:\, \mathbb{Z}_{\ge 1}^2\to \mathbb{R}_{\ge 1}$ will serve as 
bounding functions similar to $J(m,n)$.

\begin{theorem}\label{them3}
Let $a:\mathbb{Z}_{\ge 1}\to \mathbb{Z}_{\ge 1}$ be a function such that
$1\le J_1(a(n),n)\le J(a(n),n)$, $1\le J_2(a(n),n)\le J(a(n),n)$
and the conditions 
\begin{equation*}\label{44H3DB}
B_{a(n)+1} \le D_{n},\quad 4^4H^3D_{n+1} <  B_{a(n)}^3
\end{equation*} 
hold for all $n\ge n_0\ge 1$.
Assume
$\big( \overline{S}(a(n_k),n_k) \big)_{k=1}^{\infty}$ 
is an infinite sequence of non-zero solutions of equation \eqref{pohjayhtaloank} such that
\begin{equation*}\label{JJaarettomanmonta}
|x_{a(n_k),n_k}|,|y_{a(n_k),n_k}| \le H\cdot J_1(a(n_k),n_k),\quad |z_{a(n_k),n_k}|,|v_{a(n_k),n_k}| \le H\cdot J_2(a(n_k),n_k)
\end{equation*}
with an absolute constant $H\ge 1$. 
If
\begin{equation}\label{DKAPPAbound}
J_1(a(n_k),n_k)^2 \cdot J_2(a(n_k),n_k) \le D_{n_k+1}^{1-\kappa(n_k)}
\end{equation}
for all $k\in\mathbb{Z}_{\ge 1}$, then
$\alpha$ and $\beta$ satisfy the Littlewood conjecture $\ell(\alpha,\beta) = 0$.
\end{theorem}

Proof of Theorem \ref{them3}.
The proof will follow the step marks laid down in the proof of Theorem \ref{them2}.
First we note
\begin{equation*}\label{qkconstruction3}
q_k := x_{a(n_k),n_k} B_{a(n_k)} + y_{a(n_k),n_k} B_{a(n_k)+1} = z_{a(n_k),n_k} D_{n_k} + v_{a(n_k),n_k} D_{n_k+1} \ge 1
\end{equation*}
for all $n_k\ge n_0$. 
Write then $m_k=a(n_k)$, $q_k = xB_{m_k} + yB_{m_k+1} = zD_{n_k} + vD_{n_k+1}$, $J_1=J_1(a(n_k),n_k)$ and
$J_2=J_2(a(n_k),n_k)$.
Thus we may estimate
\begin{equation*}\label{}
\begin{split}
\ell_{q_k}(\alpha,\beta) 
& \le q_k\,\left(\frac{|x|}{B_{m_k+1}} + \frac{|y|}{B_{m_k+2}}\right)\,\left(\frac{|z|}{D_{n_k+1}} + \frac{|v|}{D_{n_k+2}}\right) \\
& \le 2H J_1B_{m_k+1}\,\left(\frac{HJ_1}{B_{m_k+1}} + \frac{HJ_1}{B_{m_k+2}}\right)\,\left(\frac{HJ_2}{D_{n_k+1}} + \frac{HJ_2}{D_{n_k+2}}\right) \\
& \le 8H^3 J_1B_{m_k+1}\,\frac{J_1}{B_{m_k+1}} \frac{J_2}{D_{n_k+1}}  \\
& \le 8H^3\,\frac{J_1^2J_2}{D_{n_k+1}}
  \le \frac{8H^3}{D_{n_k+1}^{\kappa(n_k)}}.
\end{split}
\end{equation*} 
Hence,
\begin{equation*}\label{liminfJJ}
\ell(\alpha,\beta) = \underset{q\to\infty}\liminf\, \ell_q(\alpha,\beta) 
\le \underset{k\to\infty}\liminf\, \frac{8H^3}{D_{n_k+1}^{\kappa(n_k)} } 
= \underset{k\to\infty}\liminf\, \frac{8H^3}{e^{\kappa(n_k)\cdot\log D_{n_k+1}} } =0.\qed
\end{equation*}

\begin{remark}\label{}
If we set $\kappa(n)=0$, then the Thue-Siegel's lemma guarantees a solution satisfying bound \eqref{DKAPPAbound}.
Therefore it will be wise to choose $\kappa(n)$ to be as small as possible.
That would improve the possibility for the existence of a solution.
\end{remark}

\subsection{The $p$-adic case}

Now we turn to the $p$-adic case which we investigate by a completely analogous manner with the Archimedean case.
Here we apply the construction, $q = xp^{m} = zD_{n} + vD_{n+1}$, for simultaneous approximations of the terms
$q|q|_p$ and $\|q\beta\|$. First we note $q|q|_p = x$. 
Therefore the $p$-adic Littlewood quantity $\ell_{p,q}(\beta) = q\,|q|_p\,\|q\beta\|$ may be written as follows
\begin{equation*}\label{}
\ell_{p,q}(\beta) 
= q\, |q|_p\,\|q\beta\| 
= xp^{m}\, |xp^{m}|_p\, \|zD_{n}\beta + vD_{n+1}\beta\|
= x\, \|zD_{n}\beta + vD_{n+1}\beta\|.
\end{equation*} 

For $(m,n)\in\mathbb{Z}^{2}_{\ge 1}$ let $\overline{S}(m,n)=(x,z,v)=(x_{m,n},z_{m,n},v_{m,n})\in\mathbb{Z}^{3}\setminus\{\overline{0}\}$
be a non-zero integer solution of the equation
\begin{equation}\label{pohjayhtalo1p}
xp^{m} = zD_{n} + vD_{n+1}.
\end{equation}
We may assume $x\ge 0$.
In this case we find a non-zero integer solution
$\overline{S}(m,n)\in\mathbb{Z}^{3}\setminus\{\overline{0}\}$ for equation \eqref{pohjayhtalo1p} such that $x\ge 1$ and
\begin{equation*}\label{}
\left\| \overline{S}(m,n) \right\|_{\infty} = \max\{|x|,|z|,|v|\} 
\le  \left\lfloor \big( p^{m} + D_{n} + D_{n+1} \big)^{1/2} \right\rfloor =: J_p(m,n).
\end{equation*}

As above for a given function $a:\mathbb{Z}_{\ge 1}\to \mathbb{Z}_{\ge 1}$ we shall consider such solutions
$\overline{S}(a(n),n)=$ $(x_{a(n),n},z_{a(n),n},v_{a(n),n})$ of the equation
\begin{equation}\label{pohjayhtalo2p}
xp^{a(n)} = zD_{n} + vD_{n+1},\quad n\in\mathbb{Z}_{\ge 1},
\end{equation}
that there exists an increasing infinite sequence $(n_k)_{k=1}^{\infty}$ satisfying
\begin{equation*}\label{qkconstruction4}
q_k := x_{a(n_k),n_k} p^{a(n_k)} = z_{a(n_k),n_k} D_{n_k} + v_{a(n_k),n_k} D_{n_k+1} \ge 1
\end{equation*}
for all $k\in\mathbb{Z}_{\ge 1}$. 
Solutions satisfying $z+v\ge 1$, $z,v\ge 0$ will be denoted by $\overline{S}^{++}(m,n)=(x,z,v)$ and called ++-solutions.

\begin{theorem}\label{thm1}
Let $a:\mathbb{Z}_{\ge 1}\to \mathbb{Z}_{\ge 1}$ be a function such that
\begin{equation*}\label{anlim0p}
\underset{n\to\infty}\lim\, \frac{ p^{a(n)} }{ D_{n} } = 0 
\end{equation*} 
and the conditions 
\begin{equation}\label{happens2pagain}
p^{a(n)} < D_{n},\quad 
3\cdot 4^2H^2 D_{n+1} < D_{n}^2
\end{equation} 
hold for all $n\ge n_0\ge 1$. 
Suppose further there exists an infinite sequence 
$\left(\overline{S}^{++}(a(n_k),n_k)\right)_{k=1}^{\infty}$ 
of ++-solutions of equation \eqref{pohjayhtalo2p} such that
\begin{equation}\label{n++aarettomanmontap}
1\le \left\| \overline{S}^{++}(a(n_k),n_k) \right\|_{\infty} \le H\cdot J_p(a(n_k),n_k)
\end{equation}
for all $k\in\mathbb{Z}_{\ge 1}$ with some absolute constant $H\ge 1$. 
Then $p$ and $\beta$ satisfy the $p$-adic Littlewood conjecture
$\ell_p(\beta) = 0$.
\end{theorem}

Proof of Theorem \ref{thm1}.
First we note Lemma \ref{LemmahJDp} with \eqref{happens2pagain} and \eqref{n++aarettomanmontap} guarantees 
\begin{equation*}\label{qkconstruction5}
q_k = x_{a(n_k),n_k} p^{a(n_k)} = z_{a(n_k),n_k} D_{n_k} + v_{a(n_k),n_k} D_{n_k+1} \ge 1
\end{equation*}
for all $n_k\ge n_0$. 
Write now $m_k=a(n_k)$, $q = xp^{m_k} =zD_{n_k} + vD_{n_k+1}$.
By the definition of $\overline{S}^{++}(m_k,n_k)$ we have $z+v\ge 1,\ z,v\ge 0$. Thereby
\begin{equation*}\label{++boundp}
(z + v)D_{n_k} \le zD_{n_k} + vD_{n_k+1} = xp^{m_k}.
\end{equation*}
By \eqref{qalphabetabound} and Lemma \ref{LemmahJDp} we get
\begin{equation*}\label{}
\|q_k\beta\|  \le \frac{|z|}{D_{n_k+1}} + \frac{|v|}{D_{n_k+2}}	
\end{equation*} 
and our assumption \eqref{n++aarettomanmontap} implies
\begin{equation*}\label{}
x \le HJ_p(m,n) \le H\left(3D_{n+1}\right)^{1/2}.
\end{equation*}
Consequently it follows
\begin{equation*}\label{}
\begin{split}
\ell_{p,q}(\beta) 
& = q\, |q|_p\,\|q\beta\| 
  = xp^{m_k}\, |xp^{m_k}|_p\, \|zD_{n_k}\beta + vD_{n_k+1}\beta\| \\
& \le x\,\left(\frac{z}{D_{n_k+1}} + \frac{v}{D_{n_k+2}}\right) 
  \le x\,\frac{z+v}{D_{n_k+1}} \\
& \le \frac{x^2p^{m_k}}{D_{n_k} D_{n_k+1}} 
  \le \frac{H^2J_p(m_k,n_k)^2p^{m_k}}{D_{n_k}D_{n_k+1}}.
\end{split}
\end{equation*}
Hence,
\begin{equation*}\label{}
\ell_{p,q}(\beta) 
\le \frac{H^2J_p(m_k,n_k)^2p^{m_k}}{D_{n_k}D_{n_k+1}}  
	\le \frac{3H^2D_{n_k+1}p^{m_k}}{D_{n_k}D_{n_k+1}} 
  =   \frac{3H^2p^{m_k}}{D_{n_k}}. 
\end{equation*}
Therefore
\begin{equation*}\label{liminfh2pankDnk}
\ell_{p}(\beta)  = \underset{q\to\infty}\liminf\, \ell_{p,q}(\beta)  
\le \underset{k\to\infty}\liminf\, \frac{3H^2p^{a(n_k)}}{D_{n_k}} = 0.\qed
\end{equation*}


\section{Effective dependence}

As noted earlier, if $B_{a(n)+1}\le D_{n+1}$, then the Thue-Siegel's lemma ensures that the equation
\begin{equation}\label{pohjayhtalo2}
xB_{a(n)} + yB_{a(n)+1}  = zD_{n} + vD_{n+1},\quad n\in\mathbb{Z}_{\ge 1},
\end{equation}
has a non-zero integer solution bounded by $\left( 4 D_{n+1} \right)^{1/3}$. 
In the sequel we will study the class of pairs of real numbers such that the solutions are infinitely often bounded by
$H D_{n+1}^{\tau}$, for $0 \le \tau <1/3$ and some constant $H\ge 1$. 
In particular, we show there are pairs of equivalent real numbers which belong to that class.
The interesting subclass with $\tau=0$ includes for example pairs of periodic type and 
pairs of palindromic continued fractions, say $\alpha=[b_0;b_{1},\ldots]$, $\beta=[d_0;d_{1},\ldots]$.
Consequently, $\|q\alpha\| \le \frac{c}{q}$ and $\|q\beta\| \le \frac{c}{q}$ hold simultaneously with 
an explicit positive constant $c$ and infinitely many positive integers $q$.
In addition, we construct examples of block-periodic continued fractions and eventually palindromic continued fractions.
Also in these cases we prove simultaneous approximation results which improve Dirichlet's theorem. 
Obviously, the Littlewood conjecture holds for such pairs.

Again we shall consider solutions
$\overline{S}(a(h),h)=(x_{a(h),h},y_{a(h),h},z_{a(h),h},v_{a(h),h})$ of equation \eqref{pohjayhtalo2}
such that there exists an increasing infinite sequence $(h_k)$ satisfying
\begin{equation*}\label{qkclassical2}
q_k := x_{a(h_k),h_k} B_{a(h_k)} + y_{a(h_k),h_k} B_{a(h_k)+1} = z_{a(h_k),h_k} D_{h_k} + v_{a(h_k),h_k} D_{h_k+1} \ge 1
\end{equation*}
for $k\in\mathbb{Z}_{\ge 1}$. 
The next theorem describes how the upper bounds of $\|q_k\alpha\|$ and $\|q_k\beta\|$ depend on the upper bound
of $\left\| \overline{S}(a(h_k),h_k) \right\|_{\infty}$.

\begin{theorem}\label{generaltau}
Let $\alpha$ and $\beta$ be real irrational numbers.
Suppose there exist a function $a:\mathbb{Z}_{\ge 1}\to \mathbb{Z}_{\ge 1}$ and 
an increasing infinite sequence of positive integers $(h_k)$ such that the solutions
$\big(\overline{S}(a(h_k),h_k)\big)_{k=1}^{\infty}$ 
of \eqref{pohjayhtalo2} satisfy
\begin{equation}\label{DleBleD}
D_{h_k} \le B_{a(h_k)+1} < D_{h_k+1},
\end{equation}
\begin{equation}\label{qk}
q_k := x_{a(h_k),h_k} B_{a(h_k)} + y_{a(h_k),h_k} B_{a(h_k)+1} = z_{a(h_k),h_k} D_{h_k} + v_{a(h_k),h_k} D_{h_k+1} \ge 1,
\end{equation}
and
\begin{equation}\label{nkaarettomanmonta}
\left\| \overline{S}(a(h_k),h_k) \right\|_{\infty} \le H \cdot D_{h_k+1}^{\tau},\quad 0\le\tau<1,
\end{equation}
for some absolute constant $H\ge 1$.  
Then the bounds
\begin{equation}\label{taubound}
\begin{split}
 \|q_k\alpha\| & \le \frac{(2H)^{2/(1+\tau)}(d_{h_k+1}+1)}{q_k^{(1-\tau)/(1+\tau)}},	\\
 \|q_k\beta\|  & \le \frac{(2H)^{2/(1+\tau)} }{q_k^{(1-\tau)/(1+\tau)}}	
\end{split}
\end{equation} 
hold for all $k\in\mathbb{Z}_{\ge 1}$. Further, there are infinitely many effectively given integers 
$q_{k}$ such that bounds in \eqref{taubound} hold.
\end{theorem}

Proof. By \eqref{qk} and \eqref{nkaarettomanmonta} we estimate
\begin{equation}\label{upperqk43tau}
q_k \le 2H D_{h_k+1}^{1+\tau}.
\end{equation}
Now the bound
$D_{h_k+1} = d_{h_k+1}D_{h_k} + D_{h_k-1} \le (d_{h_k+1}+1)D_{h_k}$ 
with \eqref{DleBleD} gives
\begin{equation*}\label{}
B_{a(h_k)+1} \ge \frac{D_{h_k+1}}{d_{h_k+1}+1}.
\end{equation*}
Combining this with 
\begin{equation*}\label{}
\|q_k\alpha\| \le \frac{|x|}{B_{a(h_k)+1}} + \frac{|y|}{B_{a(h_k)+2}} \le  \frac{2HD_{h_k+1}^{\tau}}{B_{a(h_k)+1}}
\end{equation*}
yields to
\begin{equation}\label{qkalpharepl}
\|q_k\alpha\| \le \frac{2H(d_{h_k+1}+1)}{D_{h_k+1}^{1-\tau}}.
\end{equation}
Finally, \eqref{upperqk43tau} and \eqref{qkalpharepl} give
\begin{equation*}\label{}
\|q_k\alpha\| \le \frac{(2H)^{2/(1+\tau)}(d_{h_k+1}+1)}{q_k^{(1-\tau)/(1+\tau)}}. 
\end{equation*}
For $\beta$ we have
\begin{equation}\label{tauboundDDD}
\|q_k\beta\| \le \frac{|z|}{D_{h_k+1}} + \frac{|v|}{D_{h_k+2}} \le \frac{2H}{D_{h_k+1}^{1-\tau}}
\end{equation}
instead of \eqref{qkalpharepl}.
Further, as $\tau<1$, $\beta$ is irrational and $D_{h_k+1}\to\infty$, then from \eqref{tauboundDDD} we see that
$q_k\to\infty$, too.
Hence, there are infinitely may $q_k$ such that bounds in \eqref{taubound} hold. \qed

Obviously the bounds in \eqref{taubound} are tighter than in Dirichlet's theorem when $\tau<1/3$.
At the edge, $\tau=1/3$, the situation is more subtle. Namely, by the Thue-Siegel's lemma there exists a solution 
satisfying \eqref{nkaarettomanmonta} with $H=4^{1/3}$ and $\tau=1/3$.
Therefore conditions \eqref{DleBleD} and  \eqref{qk} hold, too.
In the following we will consider this case more carefully.

Let $n\in\mathbb{Z}_{\ge 1}$ be given. Then there exists a unique $m\in\mathbb{Z}_{\ge 1}$ such that
\begin{equation}\label{13BDB}
B_{m} \le D_{n+1} < B_{m+1}.
\end{equation}
We denote this dependence by $m=a(n)$.
Note, that now for technical reasons we have chanced the roles of $B$s and $D$s.
By the Thue-Siegel's lemma there exists a solution
$\overline{S}(a(n),n)=(x_{a(n),n},y_{a(n),n},z_{a(n),n},v_{a(n),n})$ of equation \eqref{pohjayhtalo2} such that
\begin{equation}\label{qkclassical22}
q_n := x_{a(n),n} B_{a(n)} + y_{a(n),n} B_{a(n)+1} = z_{a(n),n} D_{n} + v_{a(n),n} D_{n+1} \ge 1
\end{equation}
and
\begin{equation}\label{siegel13ratkaisu}
\left\| \overline{S}(a(n),n) \right\|_{\infty} \le H\,B_{a(n)+1}^{1/3},\quad H=4^{1/3}.
\end{equation}

\begin{theorem}\label{generaltau13}
Let $\alpha$ and $\beta$ be real irrational numbers.
Then 
\begin{equation}\label{taubound13}
\begin{split}
 \|q_n\alpha\| & \le \frac{(2H)^{3/2}}{q_n^{1/2}},	\\
 \|q_n\beta\|  & \le \frac{(2H)^{3/2}(b_{a(n)+1}+1)}{q_n^{1/2}}
\end{split}
\end{equation} 
hold for all $q_n\in\mathbb{Z}_{\ge 1}$ given by \eqref{qkclassical22}. 
\end{theorem}
The proof goes in a similar manner like in Theorem \ref{generaltau} except the roles of $B$s and $D$s have changed.

\begin{corollary}\label{tau133}
Let $\alpha$ and $\beta$ be real irrational numbers, $\alpha\in\bad$ and $B=\underset{h\ge 0}\max\{b_h\}$.
Then there are infinitely many effectively given integers $q_{n}$ such that
\begin{equation}\label{taubound133}
\begin{split}
 \|q_{n}\alpha\| & \le \frac{2^{5/2}}{q_{n}^{1/2}},	\\
 \|q_{n}\beta\|  & \le \frac{2^{5/2}(B+1)}{q_{n}^{1/2}}\,.
\end{split}
\end{equation} 
\end{corollary}

Hence, we find infinitely many positive integers $q_{n}$ such that bounds in \eqref{taubound133} hold. 
The constants in the numerators are of course greater than one.
But we have an integer algorithm for computing $q_n$ for any $n\in\mathbb{Z}_{\ge 1}$.
Let $\alpha=[b_0;b_{1},\ldots]$ and $\beta=[d_0;d_{1},\ldots]$ be given continued fractions and $n\in\mathbb{Z}_{\ge 1}$. 
First compute the denominators $B_{a(n)}, B_{a(n)+1}, D_{n}$ and $D_{n+1}$ by their recurrences.
Then a solution to equation \eqref{pohjayhtalo2} may be determined in an effectively manner, say, 
by running integers $x,y,z,v$ through the interval $[-H\cdot B_{a(n)+1}^{1/3},H\cdot B_{a(n)+1}^{1/3}]$.
Therefore we may say that inequalities \eqref{taubound133} with $q_n$s given in \eqref{qkclassical22}
represent an effective version of Dirichtlet's theorem for a pair of badly approximable real numbers.

\section{Explicit constructions}

\subsection{Equivalent continued fractions}\label{secEQCF}

We shall use the notations
\begin{equation*}\label{}
\frac{A_m}{B_m}=[b_0;b_{1},\ldots,b_m],\quad \frac{C_n}{D_n}=[d_0;d_{1},\ldots,d_n],\quad \frac{E_k}{F_k}=[f_0;f_{1},\ldots,f_k]
\end{equation*}
for the convergents of $\alpha=[b_0;b_{1},\ldots]$, $\beta=[d_0;d_{1},\ldots]$ and $\gamma=[f_0;f_{1},\ldots]$, respectively.
Now we study equivalent numbers $\alpha$ and $\beta$ defined by
\begin{equation*}\label{}
\alpha = [b_0;b_1,\ldots,b_{M},f_0,f_{1},\ldots],\quad \beta = [d_0;d_1,\ldots,d_{N},f_0,f_{1},\ldots].							
\end{equation*}

Denote
\begin{equation*}\label{}
\mathcal{R}_n=
\begin{pmatrix}
b_n  & 1 \\
1    & 0
\end{pmatrix},
\qquad
\mathcal{S}_n=
\begin{pmatrix}
d_n  & 1 \\
1    & 0
\end{pmatrix},
\qquad
\mathcal{Q}_n=
\begin{pmatrix}
f_n  & 1 \\
1    & 0
\end{pmatrix}\,.  
\end{equation*}
Let $M\ge 0$ and $N\ge 0$, then we write
\begin{equation}\label{ZW}
\begin{split}
\mathcal{Z} &:= \mathcal{R}_{0}\cdots\mathcal{R}_{M} \cdot \mathcal{Q}_{0}\cdots\mathcal{Q}_{L}, \\
\mathcal{W} & := \mathcal{S}_{0}\cdots\mathcal{S}_{N} \cdot \mathcal{Q}_{0}\cdots\mathcal{Q}_{L},
\end{split} 
\end{equation}
corresponding to the convergents $[b_0;b_1,\ldots,b_{M},f_0,f_{1},\ldots,f_L]$ and 
$[d_0;d_1,\ldots,d_{N},f_0,f_{1},\ldots,f_L]$, respectively.
Otherwise $\alpha = \gamma$ or $\beta = \gamma$. In that case $\mathcal{R}_{0}\cdots\mathcal{R}_{M}=I$
or $\mathcal{S}_{0}\cdots\mathcal{S}_{N}=I$, respectively. 
Let now $M\ge 0$ and $N\ge 0$, then \eqref{ZW} corresponds to 
\begin{equation}\label{ZW2}
\begin{split}
\mathcal{Z} 
=
\begin{pmatrix}
A_{M+L+1} & A_{M+L} \\
B_{M+L+1} & B_{M+L}
\end{pmatrix}
& =
\begin{pmatrix}
A_{M} & A_{M-1} \\
B_{M} & B_{M-1}
\end{pmatrix}
\begin{pmatrix}
E_{L} & E_{L-1} \\
F_{L} & F_{L-1}
\end{pmatrix},              \\ 
\mathcal{W} 
=
\begin{pmatrix}
C_{N+L+1} & C_{N+L} \\
D_{N+L+1} & D_{N+L}
\end{pmatrix}
& =
\begin{pmatrix}
C_{N} & C_{N-1} \\
D_{N} & D_{N-1}
\end{pmatrix}
\begin{pmatrix}
E_{L} & E_{L-1} \\
F_{L} & F_{L-1}
\end{pmatrix}\,.
\end{split} 
\end{equation}

If $M\ge 0$, then we have
$\mathcal{W} = \mathcal{S}_{0}\cdots\mathcal{S}_{N} \cdot \big( \mathcal{R}_{0}\cdots\mathcal{R}_{M} \big)^{-1} \cdot \mathcal{Z}$
or equivalently
\begin{equation}\label{southeast}
\begin{split}
\begin{pmatrix}
C_{N+L+1} & C_{N+L} \\
D_{N+L+1} & D_{N+L}
\end{pmatrix}
& =  
\begin{pmatrix}
t_{M,N} & u_{M,N} \\
v_{M,N} & w_{M,N}
\end{pmatrix}
\begin{pmatrix}
A_{M+L+1} & A_{M+L} \\
B_{M+L+1} & B_{M+L}
\end{pmatrix}          \\
& =
\begin{pmatrix}
t_{M,N} A_{M+L+1} + u_{M,N} B_{M+L+1} & t_{M,N} A_{M+L} + u_{M,N} B_{M+L} \\
v_{M,N} A_{M+L+1} + w_{M,N} B_{M+L+1} & v_{M,N} A_{M+L} + w_{M,N} B_{M+L}
\end{pmatrix},
\end{split}
\end{equation}
where
\begin{equation*}\label{}
\begin{split}
\begin{pmatrix}
t_{M,N} & u_{M,N} \\
v_{M,N} & w_{M,N}
\end{pmatrix}
& = (-1)^{M-1}
\begin{pmatrix}
C_{N} & C_{N-1} \\
D_{N} & D_{N-1}
\end{pmatrix}
\begin{pmatrix}
  B_{M-1} & - A_{M-1} \\
  - B_{M} & A_{M}
  \end{pmatrix}         \\
& = (-1)^{M-1}
\begin{pmatrix}
C_{N}B_{M-1}- C_{N-1}B_{M} & - C_{N}A_{M-1} + C_{N-1}A_{M} \\
D_{N}B_{M-1}- D_{N-1}B_{M} & - D_{N}A_{M-1} + D_{N-1}A_{M}
\end{pmatrix}\,.
\end{split}
\end{equation*}
Write now
\begin{equation*}\label{writeGH}
G_{M,N} := (-1)^{M-1}(D_{N}B_{M-1}- D_{N-1}B_{M}),\quad  H_{M,N} := (-1)^{M-1}(- D_{N}A_{M-1} + D_{N-1}A_{M}).
\end{equation*}
Let $M\ge 0$. Then by the southeast corners of \eqref{southeast} we get a mixed external relation
\begin{equation}\label{EQUIVRELGEN}
D_{N+L} = G_{M,N} A_{M+L} + H_{M,N} B_{M+L}
\end{equation}
connecting $D_{N+L}$, $A_{M+L}$ and $B_{M+L}$.

Otherwise $\alpha = \gamma$. If $N\ge 0$, then
$\mathcal{W} = \mathcal{S}_{0}\cdots\mathcal{S}_{N} \cdot \mathcal{R}_{0}\cdots\mathcal{R}_{L}$
or equivalently
\begin{equation*}\label{}
\begin{split}
\begin{pmatrix}
C_{N+L+1} & C_{N+L} \\
D_{N+L+1} & D_{N+L}
\end{pmatrix}
& =
\begin{pmatrix}
C_{N} & C_{N-1} \\
D_{N} & D_{N-1}
\end{pmatrix} 
\begin{pmatrix}
A_{L} & A_{L-1} \\
B_{L} & B_{L-1}
\end{pmatrix}             \\
& =
\begin{pmatrix}
C_{N}A_{L} + C_{N-1}B_{L} & C_{N}A_{L-1} + C_{N-1}B_{L-1} \\
D_{N}A_{L} + D_{N-1}B_{L} & D_{N}A_{L-1} + D_{N-1}B_{L-1}
\end{pmatrix}\,.
\end{split}
\end{equation*}
Hence,
\begin{equation}\label{EQUIVRELVINO}
D_{N+L} = D_{N}A_{L-1} + D_{N-1}B_{L-1}.
\end{equation}
Consistent with recurrences \eqref{recurrences} one may define $A_{-1}=1$, $B_{-1}=0$, $A_{-2}=0$ and $B_{-2}=1$, 
see Continued fraction section \ref{chapter11}.
Therefore, by setting $M=-1$ we see that identity \eqref{EQUIVRELGEN} actually contains \eqref{EQUIVRELVINO}.
Further, identity \eqref{EQUIVRELVINO} at $N=0$ contains an interesting case related to $\beta = \alpha^{-1}= [0;b_0,b_1,\ldots]$ 
Namely, specializing \eqref{EQUIVRELVINO} with $D_0=1$ and $D_{-1}=0$ yields to
\begin{equation}\label{DLAL-1}
D_{L} = A_{L-1}.
\end{equation}

Note that in \eqref{EQUIVRELGEN} and \eqref{EQUIVRELVINO} the integers $M$ and $N$ are fixed per se, but $L$ is our free integer variable.

\subsection{Periodic type continued fractions}\label{pertype}

Continued fractions satisfying relation \eqref{ABrelation3} for infinitely many $h$ we call periodic type continued fractions.

Put $L=hJ-1$.
Let us suppose the elements of the matrix
\begin{equation*}\label{}
\begin{pmatrix}
E_{hJ-1} & E_{hJ-2} \\
F_{hJ-1} & F_{hJ-2}
\end{pmatrix}
\end{equation*}
connected to the continued fraction of $\gamma$, satisfy the identity
\begin{equation}\label{EFperiodic}
E_{J-2} F_{hJ-1} = F_{J-1} E_{hJ-2} 
\end{equation}
for infinitely many $h$. For example, a purely periodic continued fraction with period $J$ satisfy \eqref{EFperiodic}
for all $h\in\mathbb{Z}_{\ge 1}$, see \eqref{ABrelation3}.

Now, by \eqref{ZW2} we get
\begin{equation*}\label{}
\begin{split}
\begin{pmatrix}
E_{hJ-1} & E_{hJ-2} \\
F_{hJ-1} & F_{hJ-2}
\end{pmatrix}
& =
(-1)^{M-1}
\begin{pmatrix}
B_{M-1} & - A_{M-1} \\
- B_{M} & A_{M}
\end{pmatrix}
\begin{pmatrix}
A_{M+L+1} & A_{M+L} \\
B_{M+L+1} & B_{M+L}
\end{pmatrix}\,,
\end{split} 
\end{equation*}
which implies
\begin{equation}\label{EhJ-2FhJ-1}
\begin{split}
E_{hJ-2} & = (-1)^{M-1} (B_{M-1}A_{M+L} - A_{M-1}B_{M+L}), \\
F_{hJ-1} & = (-1)^{M-1} ( - B_{M}A_{M+L+1} + A_{M}B_{M+L+1}). 
\end{split} 
\end{equation}
Readily identities \eqref{EFperiodic} and \eqref{EhJ-2FhJ-1} give birth to an internal relation
\begin{equation}\label{tBAtAB}
F_{J-1} B_{M-1}A_{M+L} - F_{J-1} A_{M-1}B_{M+L} = - E_{J-2} B_{M}A_{M+L+1} + E_{J-2} A_{M}B_{M+L+1}. 
\end{equation}
In addition, external relation \eqref{EQUIVRELGEN} between equivalent continued fractions gives
\begin{equation*}\label{}
G_{M,N} A_{M+L} = D_{N+L} - H_{M,N} B_{M+L},
\end{equation*}
which we now insert into relation \eqref{tBAtAB} producing a purely external relation
\begin{equation}\label{PERTYPEEQUIV}
\begin{split}
& F_{J-1} B_{M-1} D_{N+L} + E_{J-2} B_{M} D_{N+L+1}
= F_{J-1} D_{N-1} B_{M+L} + E_{J-2} D_{N} B_{M+L+1} 
\end{split} 
\end{equation}
connecting $D_{N+L}$ and $D_{N+L+1}$ to $B_{M+L}$ and $B_{M+L+1}$.
Here we note that the corresponding coefficients will be constant wrt $L=hJ-1$.
Relation \eqref{PERTYPEEQUIV} produces us a common denominator $q_h$ for equivalent periodic type
continued fractions $\alpha$ and $\beta$.

\subsection{A block-periodic continued fraction}\label{blockperiodex}

Now we construct a block-periodic continued fraction which is neither periodic nor palindromic.
Remember that a block-periodic continued fraction contains arbitrary long periodic blocks. 

Put $k_n=3\cdot 4^n$, $\overset{\rightarrow} w_{0}:= abc$, and
\begin{equation}\label{wrekabcadc}
\overset{\rightarrow} w_{n+1} := \overset{\rightarrow} w_{n}\,adc\,(abc)^{k_n},\quad n=0,1,\ldots. 
\end{equation}
Further, we denote $\overset{\rightarrow} w_{\infty} := \lim \overset{\rightarrow} w_{n}$.	
The notation $t(\overset{\rightarrow} v_{n})$ will be used for the length of the word $\overset{\rightarrow} v_{n}$.
We also use $t_n:=t(\overset{\rightarrow} w_{n})$, the length of the word $\overset{\rightarrow} w_{n}$. 
Here we list some first terms
\begin{equation*}\label{}
\begin{split}											
\overset{\rightarrow} w_{0} & = abc,\quad k_0=3\cdot 1,\ t_0=3\cdot 1, \\
\overset{\rightarrow} w_{1} & = abc\, adc\, (abc)^3,\quad k_1=3\cdot 4,\ t_1=3\cdot 5, \\
\overset{\rightarrow} w_{2} & = \overset{\rightarrow} w_{1}\,adc\,(abc)^{12} 
                              = abc\, adc\, (abc)^3\, adc\, (abc)^{12},\quad k_2=3\cdot 4^2,\ t_2=3\cdot 18,\\
\overset{\rightarrow} w_{3} & = \overset{\rightarrow} w_{2}\,adc\,(abc)^{48} 
                              = abc\,adc\,(abc)^3\,adc\,(abc)^{12}\,adc\,(abc)^{48},\quad k_3=3\cdot 4^3,\ t_3=3\cdot 67.
\end{split}											
\end{equation*}
The length satisfies
\begin{equation*}\label{pituuswrekabcadc}
t_{n+1} = t_n + 3 + 3k_n,\quad t_{n} = 3(n+4^n).
\end{equation*}
In addition, we have
\begin{equation}\label{blokkipituuskn}
k_n = 3(1 + k_0 + k_1 + \ldots + k_{n-1}).
\end{equation}
Let now $a,b,c,d\in\mathbb{Z}_{\ge 1}$ be pairwise distinct.
So we may safely define an infinite continued fraction 
\begin{equation*}\label{}
\alpha:=[b_{0};b_{1},b_{2},\ldots] := [\overset{\rightarrow} w_{\infty}] =  
[a;b,c,a,d,c,a,b,c,a,b,c,a,b,c,a,d,c,a,b,c,\ldots,a,b,c,\ldots].							
\end{equation*}
Here we note that the continued fraction $\alpha$ is non-periodic. Namely, there are infinitely many $d$s whose 
successive distances are growing indefinitely.
Thus $\alpha$ is a non-quadratic irrational real number - either transcendental or algebraic of degree at least three. 
The continued fraction is also non-palindromic. If it were palindromic it should contain arbitrary long
palindromes thus containing the word $bc$. But the palindrome therefore should contain the word $cb$, which does not happen
because $a$ always follows $c$ by recurrence \eqref{wrekabcadc}.

In the following we will use the notations
\begin{equation*}\label{}
\mathcal{P}(x) =
\begin{pmatrix}
 x & 1 \\
 1 & 0
\end{pmatrix},\quad x\in\{a,b,c,d\},
\end{equation*}
\begin{equation*}\label{}
\begin{pmatrix}
A_{t_0-1} & A_{t_0-2} \\
B_{t_0-1} & B_{t_0-2}
\end{pmatrix}
=
\begin{pmatrix}
A_{2} & A_{1} \\
B_{2} & B_{1}
\end{pmatrix}
:=
\begin{pmatrix}
A_{2}(a,b,c) & A_{1}(a,b,c)  \\
B_{2}(a,b,c) & B_{1}(a,b,c) 
\end{pmatrix}
=
\mathcal{R}_{0}\mathcal{R}_{1}\mathcal{R}_{2} 
=	
\mathcal{P}(a)\mathcal{P}(b)\mathcal{P}(c),
\end{equation*}
and
\begin{equation*}\label{}
\begin{pmatrix}
\widehat A_{2} & \widehat A_{1} \\
\widehat B_{2} & \widehat B_{1}
\end{pmatrix}
:=
\begin{pmatrix}
A_{2}(a,d,c) & A_{1}(a,d,c)  \\
B_{2}(a,d,c) & B_{1}(a,d,c) 
\end{pmatrix}
=	
\mathcal{P}(a)\mathcal{P}(d)\mathcal{P}(c).
\end{equation*}
For example, the matrix product
\begin{equation*}\label{}
\begin{split}
& \mathcal{R}_{0}\mathcal{R}_{1}\mathcal{R}_{2}\cdot \mathcal{R}_{3}\mathcal{R}_{4}\mathcal{R}_{5} \cdot
\left(\mathcal{R}_{0}\mathcal{R}_{1}\mathcal{R}_{2}\right)^3
  = \mathcal{P}(a)\mathcal{P}(b)\mathcal{P}(c)\cdot\mathcal{P}(a)\mathcal{P}(d)\mathcal{P}(c) \cdot
\left(\mathcal{P}(a)\mathcal{P}(b)\mathcal{P}(c) \right)^3  \\
& =
\begin{pmatrix}
A_{2} & A_{1} \\
B_{2} & B_{1}
\end{pmatrix}
\begin{pmatrix}
\widehat A_{2} & \widehat A_{1} \\
\widehat B_{2} & \widehat B_{1}
\end{pmatrix}
\begin{pmatrix}
 A_{2} & A_{1} \\
 B_{2} & B_{1}
 \end{pmatrix}^3 
=
\begin{pmatrix}
A_{14} & A_{13} \\
B_{14} & B_{13}
 \end{pmatrix}
=
\begin{pmatrix}
A_{t_1-1} & A_{t_1-2} \\
B_{t_1-1} & B_{t_1-2}
\end{pmatrix}
\end{split}
\end{equation*}
corresponds to the convergent
\begin{equation*}\label{}
[b_{0};b_{1},b_{2},\ldots,b_{14}] = [\overset{\rightarrow} w_{1}] = [abc\, adc\, abc\, abc\, abc].							
\end{equation*}
Generally, relation \eqref{wrekabcadc} is equivalent to
\begin{equation}\label{genequiv}
\begin{split}
  \begin{pmatrix}
  A_{t_{n+1}-1}  & A_{t_{n+1}-2} \\
  B_{t_{n+1}-1}  & B_{t_{n+1}-2}
  \end{pmatrix}
 & = \mathcal{R}_{0}\cdots\mathcal{R}_{t_{n}} \cdot 
     \mathcal{P}(a)\mathcal{P}(d)\mathcal{P}(c) \cdot
     \left(\mathcal{R}_{0}\mathcal{R}_{1}\mathcal{R}_{2}\right)^{k_{n}} \\
 & = 
  \begin{pmatrix}
  A_{t_{n}-1}  & A_{t_{n}-2} \\
  B_{t_{n}-1}  & B_{t_{n}-2}
  \end{pmatrix}
\begin{pmatrix}
\widehat A_{2} & \widehat A_{1} \\
\widehat B_{2} & \widehat B_{1}
\end{pmatrix}
\begin{pmatrix}
A_{2} & A_{1} \\
B_{2} & B_{1}
\end{pmatrix}^{k_{n}}                    \\
& = 
  \begin{pmatrix}
  A_{t_{n}-1}  & A_{t_{n}-2} \\
  B_{t_{n}-1}  & B_{t_{n}-2}
  \end{pmatrix}
\begin{pmatrix}
\widehat A_{2} & \widehat A_{1} \\
\widehat B_{2} & \widehat B_{1}
\end{pmatrix}
\begin{pmatrix}
T_{3k_{n}-1} & T_{3k_{n}-2} \\ 
U_{3k_{n}-1} & U_{3k_{n}-2}
\end{pmatrix}\,,        
\end{split}
\end{equation}
where 
\begin{equation*}\label{}
\begin{pmatrix}
T_{3k_{j}-1} & T_{3k_{j}-2} \\
U_{3k_{j}-1} & U_{3k_{j}-2}
\end{pmatrix}
:=
\begin{pmatrix}
A_{2} & A_{1} \\
B_{2} & B_{1}
\end{pmatrix}^{k_{j}} 
\end{equation*}
denotes the 3-periodic block of length $3k_{j}$. 
We also use the notation
\begin{equation*}\label{}
\begin{split}
\begin{pmatrix}
Z_{k_{n}-1}  & Z_{k_{n}-2} \\
W_{k_{n}-1}  & W_{k_{n}-2}
\end{pmatrix} 
:=
&
\begin{pmatrix}
A_{2} & A_{1} \\
B_{2} & B_{1}
\end{pmatrix} 
\cdot
\begin{pmatrix}
T_{3k_{0}-1} & T_{3k_{0}-2} \\
U_{3k_{0}-1} & U_{3k_{0}-2}
\end{pmatrix}
\cdots 
\begin{pmatrix}
T_{3k_{n-1}-1} & T_{3k_{n-1}-2} \\ 
U_{3k_{n-1}-1} & U_{3k_{n-1}-2}
\end{pmatrix}\,.  
\end{split}
\end{equation*} 
By applying \eqref{11termestimates} repeatedly shows
\begin{equation}\label{ATZT}
A_{2} \prod_{j=0}^{n-1} T_{3k_{j}-1}
\le Z_{k_{n}-1}  
\le 2^{n} A_{2} \prod_{j=0}^{n-1}  T_{3k_{j}-1}.
\end{equation}
Next we observe that \eqref{blokkipituuskn} implies
\begin{equation*}\label{TV33}
\begin{split} 
\begin{pmatrix}
A_{2} & A_{1} \\
B_{2} & B_{1}
\end{pmatrix}^{k_{n}}  
& =
\left(
\begin{pmatrix}
A_{2} & A_{1} \\
B_{2} & B_{1}
\end{pmatrix}^{1}
\cdot
\begin{pmatrix}
A_{2} & A_{1} \\
B_{2} & B_{1}
\end{pmatrix}^{k_{0}} 
\cdots 
\begin{pmatrix}
A_{2} & A_{1} \\
B_{2} & B_{1}
\end{pmatrix}^{k_{n-1}}
\right)^3                          \\
& =
\left(
\begin{pmatrix}
A_{2} & A_{1} \\
B_{2} & B_{1}
\end{pmatrix}^{1}
\cdot
\begin{pmatrix}
T_{3k_{0}-1} & T_{3k_{0}-2} \\
U_{3k_{0}-1} & U_{3k_{0}-2}
\end{pmatrix}
\cdots 
\begin{pmatrix}
T_{3k_{n-1}-1} & T_{3k_{n-1}-2} \\
U_{3k_{n-1}-1} & U_{3k_{n-1}-2}
\end{pmatrix}  
\right)^3                                 \\
& =
\begin{pmatrix}
Z_{k_{n}-1}  & Z_{k_{n}-2} \\
W_{k_{n}-1}  & W_{k_{n}-2}
\end{pmatrix}^3\,. 
\end{split}
\end{equation*}
Therefore relation \eqref{wrekabcadc} is equivalent to
\begin{equation*}\label{}
\begin{split}
\begin{pmatrix}
A_{t_{n+1}-1}  & A_{t_{n+1}-2} \\
B_{t_{n+1}-1}  & B_{t_{n+1}-2}
\end{pmatrix} 
=
&
\begin{pmatrix}
A_{t_{n}-1}  & A_{t_{n}-2} \\
B_{t_{n}-1}  & B_{t_{n}-2}
\end{pmatrix}
\begin{pmatrix}
\widehat A_{2} & \widehat A_{1} \\
\widehat B_{2} & \widehat B_{1}
\end{pmatrix} 
\cdot
\begin{pmatrix}
Z_{k_{n}-1}  & Z_{k_{n}-2} \\
W_{k_{n}-1}  & W_{k_{n}-2}
\end{pmatrix}^3\,, 
\end{split}
\end{equation*}
which implies the next estimate
\begin{equation}\label{itfollowsAZ3}
A_{t_{n}-1} \widehat A_{2} Z_{k_{n}-1}^3 \le A_{t_{n+1}-1} \le 2^4 A_{t_{n}-1} \widehat A_{2} Z_{k_{n}-1}^3.
\end{equation}

On the other hand recurrence \eqref{genequiv} predicates the expansion
\begin{equation*}\label{Atnexp}
\begin{split}
\begin{pmatrix}
A_{t_{n}-1}  & A_{t_{n}-2} \\
B_{t_{n}-1}  & B_{t_{n}-2}
\end{pmatrix}
= &
\begin{pmatrix}
A_{2} & A_{1} \\
B_{2} & B_{1}
\end{pmatrix}^{1}
\cdot
\begin{pmatrix}
\widehat A_{2} & \widehat A_{1} \\
\widehat B_{2} & \widehat B_{1}
\end{pmatrix}
\begin{pmatrix}
A_{2} & A_{1} \\
B_{2} & B_{1}
\end{pmatrix}^{k_{0}} 
\cdots 
\begin{pmatrix}
\widehat A_{2} & \widehat A_{1} \\
\widehat B_{2} & \widehat B_{1}
\end{pmatrix}
\begin{pmatrix}
A_{2} & A_{1} \\
B_{2} & B_{1}
\end{pmatrix}^{k_{n-1}}                  \\                             
= &
\begin{pmatrix}
A_{2} & A_{1} \\
B_{2} & B_{1}
\end{pmatrix}^{1}
\cdot
\begin{pmatrix}
\widehat A_{2} & \widehat A_{1} \\
\widehat B_{2} & \widehat B_{1}
\end{pmatrix}
\begin{pmatrix}
T_{3k_{0}-1} & T_{3k_{0}-2} \\
U_{3k_{0}-1} & U_{3k_{0}-2}
\end{pmatrix}  
\cdots                                     \\
& 
\begin{pmatrix}
\widehat A_{2} & \widehat A_{1} \\
\widehat B_{2} & \widehat B_{1}
\end{pmatrix}
\begin{pmatrix}
T_{3k_{n-1}-1} & T_{3k_{n-1}-2} \\
U_{3k_{n-1}-1} & U_{3k_{n-1}-2}
\end{pmatrix}\,.          
\end{split}
\end{equation*}
Thus
\begin{equation}\label{AT}
A_{2}  \widehat A_{2}^n \prod_{j=0}^{n-1} T_{3k_{j}-1} \le A_{t_{n}-1} \le 2^{2n} A_{2} \widehat A_{2}^n \prod_{j=0}^{n-1}  T_{3k_{j}-1}.  
\end{equation}
By \eqref{ATZT} and \eqref{AT} we deduce
\begin{equation}\label{ZA}
\widehat A_{2} Z_{k_{n}-1} \le \widehat A_{2}^n Z_{k_{n}-1} \le 2^{n} A_{t_{n}-1}.  
\end{equation}
Further, by \eqref{itfollowsAZ3} and \eqref{ZA} we have
\begin{equation}\label{Atn+1Atn4}
A_{t_{n+1}-1} \le 2^4 A_{t_{n}-1} \widehat A_{2} Z_{k_{n}-1}^3 \le 2^{3n+4} A_{t_{n}-1}^4.  
\end{equation}
To the other direction \eqref{AT} implies
\begin{equation}\label{totheotherAtnZkn}
A_{t_{n}-1}
\le 
2^{2n} A_{2} \widehat A_{2}^n \prod_{j=0}^{n-1}  T_{3k_{j}-1}
\le 
2^{2n} \widehat A_{2}^n Z_{k_{n}-1}.
\end{equation}
Thereby  \eqref{itfollowsAZ3} and \eqref{totheotherAtnZkn} show
\begin{equation}\label{otherdir}
A_{t_{n}-1}^4
\le 
A_{t_{n}-1} 2^{6n} \widehat A_{2}^{3n} Z_{k_{n}-1}^3
\le 
2^{6n} \widehat A_{2}^{3n-1} A_{t_{n+1}-1}. 
\end{equation}

Now it is time to benefit from the block-periodicity. Namely, in the matrix product
\begin{equation}\label{benefit-block}
\begin{split}
  \begin{pmatrix}
  A_{t_{n+1}-1}  & A_{t_{n+1}-2} \\
  B_{t_{n+1}-1}  & B_{t_{n+1}-2}
  \end{pmatrix}
= 
  \begin{pmatrix}
  A_{t_{n}-1}  & A_{t_{n}-2} \\
  B_{t_{n}-1}  & B_{t_{n}-2}
  \end{pmatrix}
\begin{pmatrix}
\widehat A_{2} & \widehat A_{1} \\
\widehat B_{2} & \widehat B_{1}
\end{pmatrix}
\begin{pmatrix}
T_{3k_{n}-1} & T_{3k_{n}-2} \\
U_{3k_{n}-1} & U_{3k_{n}-2}
\end{pmatrix}       
\end{split}
\end{equation}
the last matrix represents a 3-periodic block of length $3k_{n}$.
Therefore we may apply the identity
\begin{equation*}\label{}
T_{1} U_{3k_{n}-1} = U_{2} T_{3k_{n}-2} 
\end{equation*}
implied by \eqref{ABrelation3}.
Let $d=$gcd$(T_{1},U_{2})$. Write $T_{1}=dt$, $U_{2}=ds$, where $s,t\in\mathbb{Z}_{\ge 1}$.
Then $T_{3k_{n}-2}=tZ$ and $U_{3k_{n}-1}=sZ$ for some $Z\in\mathbb{Z}_{\ge 1}$.
Consequently
\begin{equation*}\label{}
\begin{pmatrix}
T_{3k_{n}-1} & T_{3k_{n}-2} \\
U_{3k_{n}-1} & U_{3k_{n}-2}
\end{pmatrix}
=
\begin{pmatrix}
T_{3k_{n}-1} & tZ \\
sZ & U_{3k_{n}-2}
\end{pmatrix}\,,	
\end{equation*}
where $s,t$ satisfy satisfy $1\le t\le T_{1}$ and $1\le s\le U_{2}$.

Next by \eqref{benefit-block} we may express the 3-periodic block matrix as follows
\begin{equation*}\label{}
\begin{split}
&
\begin{pmatrix}
T_{3k_{n}-1} & tZ \\
sZ & U_{3k_{n}-2}
\end{pmatrix}
=
\begin{pmatrix}
\widehat A_{2} & \widehat A_{1} \\
\widehat B_{2} & \widehat B_{1}
\end{pmatrix}^{-1}
\begin{pmatrix}
A_{t_{n}-1}  & A_{t_{n}-2} \\
B_{t_{n}-1}  & B_{t_{n}-2}
\end{pmatrix}^{-1}
\begin{pmatrix}
A_{t_{n+1}-1}  & A_{t_{n+1}-2} \\
B_{t_{n+1}-1}  & B_{t_{n+1}-2}
\end{pmatrix}                               \\
& = 
(-1)^{t_{n}+3}
\begin{pmatrix}
\widehat B_{1} & - \widehat A_{1} \\
- \widehat B_{2} & \widehat A_{2}
\end{pmatrix}
\begin{pmatrix}
B_{t_{n}-2} & - A_{t_{n}-2} \\
- B_{t_{n}-1} & A_{t_{n}-1} 
\end{pmatrix}
\begin{pmatrix}
A_{t_{n+1}-1}  & A_{t_{n+1}-2} \\
B_{t_{n+1}-1}  & B_{t_{n+1}-2}
\end{pmatrix}\,.      
\end{split}
\end{equation*}
Immediately
\begin{equation*}\label{}
\begin{split}
(-1)^{t_{n}+3} tZ & = \widehat B_{1}(B_{t_{n}-2}A_{t_{n+1}-2} - A_{t_{n}-2}B_{t_{n+1}-2}) 
                              - \widehat A_{1}(- B_{t_{n}-1}A_{t_{n+1}-2} + A_{t_{n}-1}B_{t_{n+1}-2}), \\
(-1)^{t_{n}+3} sZ & = - \widehat B_{2}(B_{t_{n}-2}A_{t_{n+1}-1} - A_{t_{n}-2}B_{t_{n+1}-1}) 
                                + \widehat A_{2}(- B_{t_{n}-1}A_{t_{n+1}-1} + A_{t_{n}-1}B_{t_{n+1}-1}).
\end{split}
\end{equation*}
Let us write
\begin{equation*}\label{}
\begin{split}
&   x_n := s(\widehat B_{1}A_{t_{n}-2} + \widehat A_{1}A_{t_{n}-1}),\quad
    y_n := t(\widehat B_{2}A_{t_{n}-2} + \widehat A_{2}A_{t_{n}-1}), \\
&   z_n := s(\widehat B_{1}B_{t_{n}-2} + \widehat A_{1}B_{t_{n}-1}),\quad
    v_n := t(\widehat B_{2}B_{t_{n}-2} + \widehat A_{2}B_{t_{n}-1}).
\end{split}
\end{equation*}
Hence, we obtain an internal relation
\begin{equation}\label{ABrelation4h3}
x_n B_{t_{n+1}-2} + y_n B_{t_{n+1}-1} = z_n A_{t_{n+1}-2} + v_n A_{t_{n+1}-1}
\end{equation}
between $A$s and $B$s, with the bounds
\begin{equation*}\label{xnynznvn}
x_n, y_n, z_n, v_n \le \max\{t,s\} \widehat A_{2} A_{t_{n}-1} 
\le \max\{T_{1}, U_{2}\} \widehat A_{2} A_{t_{n}-1} .
\end{equation*}

For technical simplicity we will apply our construction \eqref{ABrelation4h3} to simultaneous approximations of 
$\alpha$ and $\beta=\alpha^{-1}$. Of course, instead we could consider any $\beta$ which is equivalent to $\alpha$, 
the block-periodic continued fraction defined above.
So let $\beta = \alpha^{-1}= [0;b_0,b_1,\ldots]$, then $D_{L} = A_{L-1}$. See Section \ref{secEQCF}, formula \eqref{DLAL-1}. 
Thereby we may define a common denominator
\begin{equation*}\label{}
q_h := x_h B_{t_{h+1}-2} + y_h B_{t_{h+1}-1} = z_h D_{t_{h+1}-1} + v_h D_{t_{h+1}}.
\end{equation*} 

\begin{theorem}\label{false periodic}
Let $\alpha$ be the block-periodic continued fraction defined above.
Then there exists an effectively determined positive constant $\gamma_4=\gamma_4(a,b,c,d)$ with respect to $h$ such that
\begin{equation*}\label{}
\begin{split}
 \|q_{h}\alpha\| & \le \frac{1}{q_h^{3/5+\delta_h}}\,, \\
 \|q_{h}\alpha^{-1}\|  & \le \frac{1}{q_h^{3/5+\delta_h}},\quad 
										   \delta_h= h\frac{\log\gamma_4}{\log q_h}\,,
\end{split}
\end{equation*} 
for all
$q_h = x_h B_{t_{h+1}-2} + y_h B_{t_{h+1}-1} = z_h D_{t_{h+1}-1} + v_h D_{t_{h+1}}$.
In addition
\begin{equation*}\label{}
\ell_{q_{h}}(\alpha,\alpha^{-1}) =\ q_{h}\,\|q_{h}\alpha\|\,\|q_{h}\alpha^{-1}\| 
\le \frac{1}{q_h^{1/5+2\delta_h}}
\end{equation*} 
for all
$h\in\mathbb{Z}_{\ge 1}$. Therefore $\alpha$ and $\alpha^{-1}$ satisfy the Littlewood conjecture
$\ell(\alpha,\alpha^{-1}) = 0$.
\end{theorem}

Proof.
Before giving estimates for $q_h$ we recall 
\begin{equation*}\label{}
\begin{split}
\phi^{t_{h}-3} & \le A_{t_{h}-1}, \\
A_{t_{h+1}-1} & \le 2^{3h+4} A_{t_{h}-1}^4 := \gamma_1^h A_{t_{h}-1}^4, \\
A_{t_{h}-1}^4 & \le 2^{6h} \widehat A_{2}^{3h-1} A_{t_{h+1}-1} := \gamma_2^h A_{t_{h+1}-1}, 
\end{split}
\end{equation*}
bounds \eqref{lowerboundAB}, \eqref{Atn+1Atn4} and \eqref{otherdir}, respectively. 
Therefore
\begin{equation}\label{3bounds}
\begin{split}
q_h & \ge A_{t_{h}-1} A_{t_{h+1}-1} \ge \phi^{t_{h}+t_{h+1}-6} \ge \phi^{4^{h+1}}, \\
q_h & \le (z_h+v_h) A_{t_{h+1}-1} \le 2\max\{T_{1}, U_{2}\} \widehat A_{2}\, A_{t_{h}-1} A_{t_{h+1}-1} 
      := \gamma_3 A_{t_{h}-1} A_{t_{h+1}-1} \\
		& \le \gamma_3 \gamma_1^h A_{t_{h}-1}^5,  \\
q_h^{3/5}	& \le (\gamma_1^h\gamma_3)^{3/5} A_{t_{h}-1}^3.
\end{split}
\end{equation} 
Now these estimates with \eqref{A=B+smallseq} yield to
\begin{equation*}\label{}
\begin{split}
 \|q_{h}\alpha\| & \le \frac{x_{h}}{B_{t_{h+1}-1}} + \frac{y_{h}}{ B_{t_{h+1}} } 
                   \le \frac{x_{h}+y_{h}}{ A_{t_{h+1}-1} } 
									 \le \frac{(b_0+1) \gamma_3 A_{t_{h}-1}}{A_{ t_{h+1}-1} }                    \\
								 & \le \frac{(a+1)\gamma_2^h\gamma_3 A_{t_{h}-1}}{A_{ t_{h}-1}^4 }
									 \le \frac{(a+1)\gamma_2^h\gamma_3 }{A_{ t_{h}-1}^3 }                      
								   \le \frac{(a+1)\gamma_2^h\gamma_3(\gamma_1^h\gamma_3)^{3/5}}{q_h^{3/5}}   \\
								 & \le \frac{\gamma_4^h}{q_h^{3/5}}                                            
								  	=  \frac{1}{q_h^{3/5+\delta_h}},\quad 
										   \delta_h= h\frac{\log\gamma_4}{\log q_h}\,,																				\\  
\end{split}
\end{equation*} 
where $\gamma_4=\gamma_4(a,b,c,d)$ is constant with respect to $h$ satisfying
$(a+1)\gamma_2^h\gamma_3(\gamma_1^h\gamma_3)^{3/5} \le \gamma_4^h$ for $h$ big enough.
Obviously $\gamma_4(a,b,c,d)$ may be given effectively.
Similarly holds an estimate
\begin{equation*}\label{}
\|q_{h} \beta\|  \le \frac{z_h}{ D_{t_{h+1}} } + \frac{v_h}{D_{t_{h+1}+1}}	
                 \le \frac{z_h+v_h}{A_{t_{h+1}-1}}
								 \le \frac{1}{q_h^{3/5+\delta_h}}\,.																			
\end{equation*}
The first bound in \eqref{3bounds} implies
\begin{equation*}\label{}
\delta_h = h\frac{\log\gamma_4}{\log q_h} \le \frac{h}{4^{h+1}}\frac{\log\gamma_4}{\log\phi}.
\end{equation*} 
Thus $\delta_h\ \to\ 0$ as  $h\ \to\ \infty$. 
Therefore, from
\begin{equation*}\label{}
\ell_{q_{h}}(\alpha,\beta) =\ q_{h}\,\|q_{h}\alpha\|\,\|q_{h}\beta\| 
\le \frac{1}{q_h^{1/5+2\delta_h}}
\end{equation*} 
follows
$\ell(\alpha,\beta) = 0$.\qed

\subsection{An eventually palindromic continued fraction}\label{eventuallypalindex}

Remember that a continued fraction is called purely palindromic (palindromic), if it has arbitrary long convergents which are palindromes.
Further, a continued fraction is called eventually palindromic, if it contains arbitrary long palindromes but
the continued fraction itself is not necessarily purely palindromic.

Next we like to show how to construct examples which are not necessarily purely palindromic but contain arbitrary long palindromes. 
For that we may apply recurrences like 
\begin{equation*}\label{wrekpp}
\overset{\rightarrow} w_{n+1} := 
\overset{\rightarrow} w_{n} \overset{\rightarrow} p_{n} \overset{\leftarrow}{p}_{n},\quad n=1,2,\ldots. 
\end{equation*}
in a monoid $\langle a_1,a_2,\ldots \rangle$ of words starting from an initial value 
$\overset{\rightarrow} w_{1}\in \langle a_1,a_2,\ldots \rangle$.
In addition in every step we select our favorite word $\overset{\rightarrow} p_{n}\in \langle a_1,a_2,\ldots \rangle$.
Further, define $\overset{\rightarrow} w_{\infty}:=\lim\overset{\rightarrow} w_{n}$.
Note that even $\overset{\rightarrow} p_{n} \overset{\leftarrow}{p}_{n}$ is a palindrome but 
$\overset{\rightarrow} w_{\infty}$ is not necessarily purely palindromic.

As an example we choose $\overset{\rightarrow} w_{1}:= ab\in\langle a,b \ldots \rangle$ and
$\overset{\rightarrow} p_{n} = \overset{\rightarrow} w_{n} \overset{\rightarrow} w_{n}$.
Thus the recurrence looks like
\begin{equation}\label{wrek5}
\overset{\rightarrow} w_{n+1} 
:= \overset{\rightarrow} w_{n}  \overset{\rightarrow} w_{n} \overset{\rightarrow} w_{n}  \overset{\leftarrow} w_{n} \overset{\leftarrow} w_{n}
 = \overset{\rightarrow} w_{n} \overset{\rightarrow} w_{n}^2 \overset{\leftarrow} w_{n}^2 
\end{equation}
with the initial value $\overset{\rightarrow} w_{1}:= ab$. 
Some first terms are given by
$\overset{\rightarrow} w_{1} = ab$, 
$\overset{\rightarrow} w_{2} = abababbaba$,
$\overset{\rightarrow} w_{3} = abababbaba\, abababbaba$ $abababbaba\, ababbababa\, ababbababa$. 
Note that the word $\overset{\rightarrow} w_{n}^2 \overset{\leftarrow} w_{n}^2$ is a palindrome. 
Let now $a,b\in\mathbb{Z}_{\ge 1}$, $a\ne b$. So we may define an infinite eventually palindromic continued fraction 
\begin{equation}\label{abababbaba...}
\alpha:=[r_{0};r_{1},r_{2},\ldots] := [\overset{\rightarrow} w_{\infty}] = [a;b,a,b,a,b,b,a,b,a,\ldots].							
\end{equation}

The continued fraction \eqref{abababbaba...} is a non-quadratic number because the word 
$\overset{\rightarrow} w_{\infty}$ is non-periodic, see Lemma \ref{5wnonperiodic}.
We want to focus on our main theme, explicit simultaneous approximations and construction of common denominators.
Therefore we postpone proof of Lemma \ref{5wnonperiodic} until the main result, Theorem \ref{false palindromic}, of this section.

Turning to matrix representations we see relation \eqref{wrek5} is equivalent to 
\begin{equation}\label{matrix5h+4}
\begin{split}
  \begin{pmatrix}
  A_{5h+4}  & A_{5h+3} \\
  B_{5h+4}  & B_{5h+3}
  \end{pmatrix}
 & = \big(\mathcal{R}_{0}\cdots\mathcal{R}_{h}\big)^3  \cdot \big(\mathcal{R}_{h}\cdots\mathcal{R}_{0}\big)^2 \\
 & = 
  \begin{pmatrix}
  A_{h}  & A_{h-1} \\
  B_{h}  & B_{h-1}
  \end{pmatrix}^3
  \begin{pmatrix}
  A_{h}   & B_{h} \\
  A_{h-1} & B_{h-1}
  \end{pmatrix}^2	          \\
& = 
  \begin{pmatrix}
  A_{h}  & A_{h-1} \\
  B_{h}  & B_{h-1}
  \end{pmatrix}
	\begin{pmatrix}
  A_{2h+1}  & A_{2h} \\
  B_{2h+1}  & B_{2h}
  \end{pmatrix}
 \begin{pmatrix}
  A_{2h+1}  & B_{2h+1} \\
  A_{2h}  & B_{2h}
  \end{pmatrix}         \\
& =: 
  \begin{pmatrix}
  A_{h}  & A_{h-1} \\
  B_{h}  & B_{h-1}
  \end{pmatrix}
 \begin{pmatrix}
  Y & Z \\
  Z & W
  \end{pmatrix}\,,	          
\end{split}
\end{equation}
where we have applied palindrome property \eqref{palindromematrix}.
Readily we get
\begin{equation*}\label{}
\begin{split}
(-1)^{h-1}
\begin{pmatrix}
  B_{h-1} & - A_{h-1} \\
  - B_{h} & A_{h}
  \end{pmatrix}
  \begin{pmatrix}
  A_{5h+4}  & A_{5h+3} \\
  B_{5h+4}  & B_{5h+3}
  \end{pmatrix}
	=
 \begin{pmatrix}
  Y & Z \\
  Z & W
  \end{pmatrix}	          
\end{split}
\end{equation*}
and further
\begin{equation*}\label{}
\begin{split}
(-1)^{h-1}\big(B_{h-1} A_{5h+3} - A_{h-1} B_{5h+3}\big) &	= Z, \\
(-1)^{h-1}\big(- B_{h} A_{5h+4} + A_{h} B_{5h+4}\big) &	= Z.
\end{split}
\end{equation*}
Hence, we obtain an internal relation
\begin{equation*}\label{esiq5h}
A_{h} B_{5h+4} + A_{h-1} B_{5h+3} = B_{h} A_{5h+4} + B_{h-1} A_{5h+3},
\end{equation*}
which allows us to define a Janus face common denominator 
\begin{equation*}\label{}
q_{h} := A_{h} B_{5h+4} + A_{h-1} B_{5h+3} = B_{h} D_{5h+5} + B_{h-1} D_{5h+4},\quad h\in\mathbb{Z}_{\ge 1},
\end{equation*}
for $\alpha$ and $\beta=\alpha^{-1}$.

Again we will apply our construction to simultaneous approximation of $\alpha$ and $\beta=\alpha^{-1}$.
So, the proof of the next theorem goes in a similar fashion like in Theorem \ref{false periodic}.

\begin{theorem}\label{false palindromic}
Let $\alpha$ be the eventually palindromic continued fraction defined above.
Then  
\begin{equation*}\label{}
\begin{split}
 \|q_{h}\alpha\| & \le 2(a+1) \left(\frac{32}{a}\right)^{2/3} \frac{1}{q_h^{2/3}}\,,	\\
 \|q_{h}\alpha^{-1}\|  & \le \frac{2}{a} \left(\frac{32}{a}\right)^{2/3} \frac{1}{q_h^{2/3}}
\end{split}
\end{equation*} 
for all
$q_{h} := A_{h} B_{5h+4} + A_{h-1} B_{5h+3} = B_{h} D_{5h+5} + B_{h-1} D_{5h+4}, h\in\mathbb{Z}_{\ge 1}$.  
Further
\begin{equation*}\label{}
\ell_{q_{h}}(\alpha,\alpha^{-1}) =\ q_{h}\,\|q_{h}\alpha\|\,\|q_{h}\alpha^{-1}\| 
\le \frac{4(a+1) }{a} \left(\frac{32}{a}\right)^{4/3} \frac{1}{q_{h}^{1/3}}
\end{equation*} 
for all
$q_{h}$, $h\in\mathbb{Z}_{\ge 1}$. Therefore $\alpha$ and $\alpha^{-1}$ satisfy the Littlewood conjecture
$\ell(\alpha,\alpha^{-1}) = 0$.
\end{theorem}

Proof. Let us start by noting that \eqref{A=B+smallseq} with $b_0=a$ gives 
$aB_{k} < A_{k} < (a+1)B_{k}$ for $k\in\mathbb{Z}_{\ge 1}$.
Thereby matrix product \eqref{matrix5h+4} implies the estimates
\begin{equation*}\label{AhA5h+4}
A_h^5 \le A_{5h+4} \le 16A_h^5,\quad 
q_h   \le \frac{2}{a}A_hA_{5h+4} \le \frac{32}{a}A_h^6.
\end{equation*} 
Thus
\begin{equation*}\label{}
\begin{split}
 \|q_{h}\alpha\| & \le \frac{A_{h}}{B_{5h+5}} + \frac{A_{h-1}}{B_{5h+4}}
                   \le \frac{2(a+1)A_{h}}{A_{5h+4}} \le \frac{2(a+1)}{A_{h}^4}
                   \le 2(a+1) \left(\frac{32}{a}\right)^{2/3} \frac{1}{q_h^{2/3}},\\
 \|q_{h}\beta\|  & \le \frac{B_{h}}{D_{5h+6}} + \frac{B_{h-1}}{D_{5h+5}} 
                   \le \frac{2B_{h}}{D_{5h+5}} \le \frac{2A_{h}}{aA_{5h+5}}  
								   \le \frac{2}{a} \frac{1}{A_{h}^4}
								   \le \frac{2}{a} \left(\frac{32}{a}\right)^{2/3} \frac{1}{q_h^{2/3}}
\end{split}
\end{equation*} 
for $h\in\mathbb{Z}_{\ge 1}$.
Finally
\begin{equation*}\label{}
\ell_{q_{h}}(\alpha,\beta) =\ q_{h}\,\|q_{h}\alpha\|\,\|q_{h}\beta\| 
\le \frac{4(a+1) }{a} \left(\frac{32}{a}\right)^{4/3} \frac{1}{q_{h}^{1/3}},
\end{equation*} 
which shows
$\ell(\alpha,\alpha^{-1}) = 0$.\qed

Before showing the word $\overset{\rightarrow} w_{\infty}$ is non-periodic we will present some fundamental properties
which are applications of recurrence \eqref{wrek5}.
The length of $\overset{\rightarrow} w_{n}$ satisfies 
$t(\overset{\rightarrow} w_{n}) = 2\cdot 5^{n-1}$.
For $k\in\mathbb{Z}_{\ge 2}$, we have the following representations
\begin{equation}\label{wrekk}
\overset{\rightarrow} w_{\infty} 
= \prod_{j=1}^{\infty} \overset{\rightarrow} w_{k}\overset{\rightarrow} w_{k}\overset{\leftrightarrow} w_{k,j}
	                       \overset{\leftarrow} w_{k}\overset{\leftarrow} w_{k}
= \prod_{j=1}^{\infty} \overset{\leftrightarrow} w_{k,j}, 
\end{equation}
where $\overset{\leftrightarrow} w_{k,j}\in \{\overset{\rightarrow} w_{k},\overset{\leftarrow} w_{k}\}$.
We may further write 
$\overset{\rightarrow}w_{\infty} = \overset{\rightarrow}w_{h} K_h,\ 
K_h = \overset{\rightarrow}w_{h}\overset{\rightarrow}w_{h}
      \overset{\leftarrow}w_{h}\overset{\leftarrow}w_{h}\, \overset{\rightarrow}w_{h}\cdots$
for $h\in\mathbb{Z}_{\ge 1}$.
The middle factor of $\overset{\rightarrow} w_{n}$ is $ab$ and the middle factor of $\overset{\leftarrow} w_{n}$ is $ba$. 
Therefore
$\overset{\rightarrow} w_{n} \ne \overset{\leftarrow} w_{n}$
for all $n\in\mathbb{Z}_{\ge 1}$.
In the sequel the relation $v\sqsubseteq w$ tells $v$ is a factor (substring) of the word $w$. 
The empty word is denoted by $\epsilon$.

\begin{lemma}\label{5wnonperiodic}
The word $\overset{\rightarrow} w_{\infty}$ is non-periodic.
\end{lemma}\label{}
Proof.
Suppose $P'\sqsubseteq\overset{\rightarrow}w_{\infty}$ is a period of length $J$. 
Then there exists a $k$ such that 
$P'\sqsubseteq \overset{\rightarrow} w_{k}= \overset{\rightarrow} v_{1}P' \overset{\rightarrow} v_{2}$, where
$\overset{\rightarrow} w_{k}\sqsubseteq \overset{\rightarrow} w_{\infty}$ is a prefix of $\overset{\rightarrow} w_{\infty}$. 
By looking at
\begin{equation}\label{secondblock}
\overset{\rightarrow} w_{k+1} 
= \overset{\rightarrow} w_{k}  \overset{\rightarrow} w_{k} \overset{\rightarrow} w_{k}  \overset{\leftarrow} w_{k} \overset{\leftarrow} w_{k}
= \overset{\rightarrow} v_{1} P'P'P'\cdots
\end{equation}
we see that there exists a period $P$ of length $J$, which is a prefix of 
$\overset{\rightarrow} w_{k} \overset{\rightarrow} w_{k} \overset{\leftarrow} w_{k} \overset{\leftarrow} w_{k}$.
Hence, the period candidate $P$ will be a prefix in the second block $\overset{\rightarrow} w_{k}$ determined by
representation \eqref{secondblock}.
It follows the period
$P \sqsubseteq \overset{\rightarrow} w_{k}\sqsubseteq K_k$
is a prefix of $K_k$ and more precisely a prefix of the first $\overset{\rightarrow} w_{k}$ in $K_k$.

Let $2\le h\le k$. The first prefix $\overset{\rightarrow} w_{h}$ of $K_h$ is not a period.
Suppose on the contrary that $\overset{\rightarrow}w_{h}$ is a period $P$. Then
$K_h = PPPPP\cdots = \overset{\rightarrow} w_{h}\overset{\rightarrow}w_{h}
                    \overset{\leftarrow} w_{h}
	                  \overset{\leftarrow} w_{h}\, \overset{\rightarrow}w_{h}\cdots$
implies $\overset{\rightarrow}w_{h}=\overset{\leftarrow}w_{h}$. A contradiction.

Suppose $P\sqsubseteq\overset{\rightarrow}w_{2}$ and $J<10$. 
We have 
$K_2 = \overset{\rightarrow}w_{2}\overset{\rightarrow}w_{2}\overset{\leftarrow}w_{2}\overset{\leftarrow} w_{2}\cdots
= abababbaba\,abababbaba \cdots.$
Going trough all possibilities
$a$, $ab$,\ldots, $abababbab$, we see that $J\ge 10$.
If $J=10$, then $\overset{\rightarrow} w_{2}$ would be a period. A contradiction.
Hence, there is no period of length $J\le 10$. In particular $\overset{\rightarrow} w_{2}$ is not a period. 
Therefore, if $P\sqsubseteq\overset{\rightarrow} w_{k}$, then $k\ge 3$.

Next we prove  $10|J$. First we show $5|J$.
Suppose on the contrary $5\not|J$. 
For technical reason we will use the following indexing for the partial quotients in
$\alpha:=[s_{1};s_{2},s_{3},\ldots] := [\overset{\rightarrow} w_{\infty}] = [a;b,a,b,a,b,b,a,b,a,\ldots]$.
We start by noting
$\overset{\rightarrow} w_{\infty} 
=  s_1s_2s_3\cdots
=  \prod_{j=1}^{\infty} \overset{\rightarrow} w_{2}\overset{\rightarrow} w_{2}\overset{\leftrightarrow} w_{2,j}
	                       \overset{\leftarrow} w_{2}\overset{\leftarrow} w_{2}$.
The period $P$ is a prefix of the second $\overset{\rightarrow} w_{k}$, $k\ge 3$, in product \eqref{wrekk}. Thus 
$P  =s_{10h+1}s_{10h+2}\cdots s_{10h+J} 
    = \overset{\rightarrow} w_{2}\overset{\rightarrow} w_{2}\overset{\rightarrow} w_{2}\overset{\leftarrow} w_{2}\overset{\leftarrow} w_{2}\cdots
    = abababbaba\,aba \cdots$
for some $h\in \mathbb{Z}_{\ge 1}$.												
Therefore
\begin{equation*}\label{}
s_{10h+1}s_{10h+2}\cdots s_{10h+10} =
\begin{cases}											
abababbaba\quad\text{or} \\	
ababbababa
\end{cases}											
\end{equation*}
for $h\in \mathbb{Z}_{\ge 1}$. 
In particular, $s_{10h+3}=s_{10h+8}=a$ and $s_{10h+4}=s_{10h+9}=b$ for all $h\in \mathbb{Z}_{\ge 1}$.
Because $5\not|J$ there exist $m,m'\in \mathbb{Z}_{\ge 1}$ such that
\begin{equation*}\label{}
\begin{cases}											
mJ \equiv 3=a \pmod{5}; \\	
m'J\equiv 4=b \pmod{5}.
\end{cases}											
\end{equation*}
It follows $10mh+mJ=10H+3$ or $10mh+mJ=10H+8$ and $10m'h+m'J=10H'+4$ or $10m'h+m'J=10H'+8$.
Hence, the last digits of the periodic blocks $P^{m}$ and $P^{m'}$ will be
$s_{10mh+mJ} = a$ and $s_{10mh+m'J} = b$,
respectively. A contradiction with the fact that $P^{m}$ and $P^{m'}$ should have identical last digits.\\
Suppose then $J=(2l+1)5$. Then in the middle of the word $PP$ we have either $babb$ or $bbab$
while in the middle of the word $P^{2}P^{2}$ we always have $baab$. Therefore $10|J$.\qed\\    

Again we apply the presentation
$\overset{\rightarrow}w_{\infty} = \overset{\rightarrow}w_{k}\, K_k, 
K_k = \overset{\rightarrow}w_{k}\overset{\rightarrow}w_{k}
\overset{\leftarrow} w_{k}\overset{\leftarrow}w_{k}\, \overset{\rightarrow}w_{k}\cdots.$
Because $J=10H$ for some $H\in \mathbb{Z}_{\ge 1}$ it follows that the period $P$ has a representation as a product
$P = \prod_{j=1}^{H} \overset{\leftrightarrow} w_{2,j}\sqsubseteq \overset{\rightarrow} w_{k}\sqsubseteq K_k$
and $P$ is a prefix of $K_k$ and more precisely a prefix of the first $\overset{\rightarrow} w_{k}$ in $K_k$.
We may assume that $k\ge 3$ is the smallest possible satisfying the above inclusion.
Let us list these crucial facts to be used in the sequel
\begin{equation*}\label{}
\begin{split}
\overset{\rightarrow} w_{\infty} = & \overset{\rightarrow} w_{k}\,K_k,\quad  
\overset{\rightarrow} w_{k-1}\sqsubseteq P \sqsubseteq \overset{\rightarrow} w_{k},\quad
\overset{\rightarrow} w_{k-1}\ne P,                                                       \\ 
K_k = & PPPPPPP\cdots                                                                         
    =   \overset{\rightarrow} w_{k}\overset{\rightarrow} w_{k}
		    \overset{\leftarrow} w_{k}
	      \overset{\leftarrow} w_{k}\, \overset{\rightarrow} w_{k}\cdots
		=   \overset{\rightarrow} w_{j}\overset{\rightarrow} w_{j}\overset{\rightarrow} w_{j}
	      \overset{\leftarrow} w_{j} \overset{\leftarrow} w_{j}\, 
			  \overset{\rightarrow} w_{j}\overset{\rightarrow} w_{j}\overset{\rightarrow} w_{j}
	      \overset{\leftarrow} w_{j} \overset{\leftarrow} w_{j}\, 	
				\cdots   																																										
\end{split}			
\end{equation*}
for $2\le j\le k-1$.

As seen above $\overset{\rightarrow} w_{2}$ is not a period.
Therefore we write the period as a product
\begin{equation*}\label{}
\begin{split}
P =    & \prod_{j=1}^{H_3} \overset{\leftrightarrow} w_{3,j}\,E_3
  =      \overset{\rightarrow}w_{3}\overset{\rightarrow}w_{3}\overset{\rightarrow}w_{3}\overset{\leftarrow} w_{3}\overset{\leftarrow} w_{3} 
         \cdots 
	       \overset{\leftrightarrow} w_{3} \overset{\leftrightarrow} w_{3}\,E_3,																						
\end{split}			
\end{equation*}
where $E_3\sqsubseteq\overset{\leftrightarrow} w_{3}$, 
$E_3\in 	\{ \epsilon,
             \overset{\rightarrow} w_{2},
             \overset{\rightarrow} w_{2}\overset{\rightarrow} w_{2},
     		  	\overset{\rightarrow} w_{2}\overset{\rightarrow} w_{2}\overset{\leftrightarrow} w_{2},
     				\overset{\rightarrow} w_{2}\overset{\rightarrow} w_{2}\overset{\leftrightarrow} w_{2}\overset{\leftarrow} w_{2} \}$,
a prefix of $\overset{\leftrightarrow}w_{3}$, is a suffix of the first period $P$.
The corresponding non-empty suffixes 
$\{\overset{\rightarrow} w_{2}\overset{\leftrightarrow} w_{2}\overset{\leftarrow} w_{2}\overset{\leftarrow} w_{2},
  \overset{\leftrightarrow} w_{2}\overset{\leftarrow} w_{2}\overset{\leftarrow} w_{2}, 
  \overset{\leftarrow} w_{2}\overset{\leftarrow} w_{2},											
 	\overset{\leftarrow} w_{2} \}$
of $\overset{\leftrightarrow}w_{3}$ determine the prefixes of the next period, which should start as the first period
$P = \overset{\rightarrow} w_{2}\overset{\rightarrow} w_{2}\overset{\rightarrow} w_{2}
	   \overset{\leftarrow} w_{2}\overset{\leftarrow} w_{2}\cdots$.
Obviously there are no matches. 
In the empty word case we have
\begin{equation*}\label{}
\begin{split}
P = & \prod_{j=1}^{H_3} \overset{\leftrightarrow} w_{3,j}
  = \prod_{j=1}^{H_4} \overset{\leftrightarrow} w_{4,j}\,E_4,\quad
E_4\in \{ \epsilon,
          \overset{\rightarrow} w_{3},
          \overset{\rightarrow} w_{3}\overset{\rightarrow} w_{3},
     		  \overset{\rightarrow} w_{3}\overset{\rightarrow} w_{3}\overset{\leftrightarrow} w_{3},
     			\overset{\rightarrow} w_{3}\overset{\rightarrow} w_{3}\overset{\leftrightarrow} w_{3}\overset{\leftarrow} w_{3} \},																							
\end{split}			
\end{equation*}
where $E_4\sqsubseteq\overset{\leftrightarrow} w_{4}$ is a suffix of the first period $P$.

Recursively we reach the representation
\begin{equation*}\label{}
P = \prod_{j=1}^{H_{k-1}} \overset{\leftrightarrow} w_{k-1,j}\,E_{k-1}, 
\end{equation*}
where $E_{k-1}\sqsubseteq\overset{\leftrightarrow} w_{k-1}$ is a suffix of the first period $P$.
By the facts 
$\overset{\rightarrow} w_{k-1}\sqsubseteq P \sqsubseteq \overset{\rightarrow} w_{k}$, 
$\overset{\rightarrow} w_{k-1}\ne P$ and
$\overset{\rightarrow} w_{k} 
= \overset{\rightarrow}w_{k-1}\overset{\rightarrow}w_{k-1}\overset{\rightarrow}w_{k-1}\overset{\leftarrow}w_{k-1}\overset{\leftarrow}w_{k-1}$
we see that
\begin{equation*}\label{}
E_{k-1}\in \{ \epsilon,
                 \overset{\rightarrow} w_{k-2},
                 \overset{\rightarrow} w_{k-2}\overset{\rightarrow} w_{k-2},
     		  	     \overset{\rightarrow} w_{k-2}\overset{\rightarrow} w_{k-2}\overset{\leftrightarrow} w_{k-2},
     				     \overset{\rightarrow} w_{k-2}\overset{\rightarrow} w_{k-2}\overset{\leftrightarrow} w_{k-2}\overset{\leftarrow} w_{k-2} \}.
\end{equation*} 
The non-empty suffixes $E_{k-1}$ determine the prefixes
\begin{equation*}\label{}
\begin{split}
\{\overset{\rightarrow} w_{k-2}\overset{\leftrightarrow} w_{k-2}\overset{\leftarrow} w_{k-2}\overset{\leftarrow} w_{k-2},
  \overset{\leftrightarrow} w_{k-2}\overset{\leftarrow} w_{k-2}\overset{\leftarrow} w_{k-2}, 
  \overset{\leftarrow} w_{k-2}\overset{\leftarrow} w_{k-2},											
 	\overset{\leftarrow} w_{k-2} \}																							
\end{split}			
\end{equation*}
of the next period, which should start as the first period
$P = \overset{\rightarrow}w_{k-2}\overset{\rightarrow}w_{k-2}\overset{\rightarrow}w_{k-2}
	   \overset{\leftarrow}w_{k-2}\overset{\leftarrow}w_{k-2}\cdots$.
There are no matches.
In the case of the empty word we have the following possibilities
\begin{equation*}\label{}						
P\in \{ \overset{\rightarrow}w_{k-1}\overset{\rightarrow}w_{k-1},
        \overset{\rightarrow}w_{k-1}\overset{\rightarrow}w_{k-1}\overset{\rightarrow} w_{k-1},
        \overset{\rightarrow}w_{k-1}\overset{\rightarrow}w_{k-1}\overset{\rightarrow} w_{k-1} \overset{\leftarrow} w_{k-1},
	      \overset{\rightarrow} w_{k} \}
\end{equation*} 
for the first period $P\sqsubseteq\overset{\rightarrow}w_{k}$, which is a prefix of 
\begin{equation}\label{}						
K_{k} = PPPP\cdots = \overset{\rightarrow}w_{k-1}\overset{\rightarrow}w_{k-1}\overset{\rightarrow}w_{k-1}
                     \overset{\leftarrow}w_{k-1}\overset{\leftarrow}w_{k-1}\,
										 \overset{\rightarrow}w_{k-1}\overset{\rightarrow}w_{k-1}\overset{\rightarrow}w_{k-1}
                     \overset{\leftarrow}w_{k-1}\overset{\leftarrow}w_{k-1}\,\cdots.
\end{equation}
Again a contradiction.\qed

\emph{Acknowledgements.}
The author is indebted to Topi T\"orm\"a of his valuable remarks.
Further thanks to Camilla Hollanti for her hospitality.
The research was partially supported by the Academy of Finland grant \#351271.

\end{document}